\documentclass[11pt,twoside]{article}
\usepackage{latexsym,amssymb,amsmath}
\usepackage{amsfonts}
\RequirePackage{graphics,psfrag}
\def\abs#1{\left \vert #1 \right \vert}

\usepackage{multirow}
\def\abs#1{\left \vert #1 \right \vert}
\def\R{{\bf R}} 
\def\Q{{\bf Q}} 
\def\SS{{\bf S}} 

\def\pn{\medskip\par\noindent}
\def\Frac#1#2{{\displaystyle{{#1}\over{#2}}}}
\def\[#1\]{\begin{eqnarray}#1\end{eqnarray}}
\def\$#1\${\begin{eqnarray}#1\end{eqnarray}}
\def\eps{\varepsilon}

\def\Mod#1{\,(\hbox{\rm mod}\,#1)}
\def\Res{\hbox{\rm Res}\,}
\def\bi{\begin{itemize}}
\def\ei{\end{itemize}}
\def\bn{\begin{enumerate}}
\def\en{\end{enumerate}}
\def\cC{{\cal C}}

\def\cA{{\cal A}}
\def\cB{{\cal B}}
\def\cZ{{\cal Z}}

\def\cK{{\cal K}}
\def\sign#1{{\rm sign}\,\bigl( #1 \bigr)}

\def\abs#1{\left \vert #1 \right \vert}
\def\Frac#1#2{{\displaystyle{\frac{#1}{#2}}}}
\def\phi{\varphi}
\catcode`\@=11
\def\@opargbegintheorem#1#2#3{\par\addvspace{6pt plus3pt minus2pt}%
    \def\@tempa{#3}%
    \noindent{\bf #1 #2 \ifx\@tempa\empty\unskip\else\unskip\ (#3).\fi\hskip.5em}\csname#1font\endcsname\ignorespaces
\ignorespaces}

\def\@endtheorem{\par\addvspace{6pt plus3pt minus2pt}}
\def\@begintheorem#1#2#3{\par\addvspace{8pt plus3pt minus2pt}%
              \noindent{\csname#1headfont\endcsname#1\ \ignorespaces#3 #2.}%
              \csname#1font\endcsname\hskip6pt\ignorespaces}
\def\@endtheorem{\par\addvspace{8pt plus3pt minus2pt}\@endparenv}

\catcode`\@=12

\newtheorem{theorem}{Theorem}[section]
\newtheorem{thm*}{Theorem}
\newtheorem{cor}[theorem]{Corollary}
\newtheorem{lemma}[theorem]{Lemma}
\newtheorem{proposition}[theorem]{Proposition}

\newtheorem{remark}[theorem]{Remark}

\newtheorem{algorithm}[theorem]{Algorithm}
\newcommand{\Pf}{{\em Proof}. }

\newcommand{\EPf}{\hbox{}\hfill$\Box$\vspace{.5cm}}
\begin{document}
\title{The first rational Chebyshev knots\footnote{Expanded version of a conference
presented at Mega 09 (Barcelona, June 2009).}}
\author{P. -V. Koseleff, D. Pecker, F. Rouillier}
\maketitle
\begin{abstract}
A Chebyshev knot $\cC(a,b,c,\phi)$ is a knot which has a
parametrization of the form $ x(t)=T_a(t); \  y(t)=T_b(t) ; \
z(t)= T_c(t + \phi), $ where $a,b,c$ are integers, $T_n(t)$ is the
Chebyshev polynomial of degree $n$ and $\phi \in \R.$ We show that
any two-bridge knot is a Chebyshev knot with $a=3$ and also with
$a=4$. For every $a,b,c$ integers ($a=3, 4$ and $a$, $b$ coprime), we
describe an algorithm that gives all Chebyshev knots $\cC(a,b,c,\phi)$.
We deduce a list of minimal Chebyshev representations of two-bridge
knots with small crossing number.
\end{abstract}
\pn {\bf keywords:} Polynomial curves; two-bridge knots;
Chebyshev curves; real roots isolation;
computer algebra; algorithms
\section{Introduction}
A Chebyshev knot $\cC(a,b,c,\phi)$ is a knot which has a parametrization of the
form $ x(t)=T_a(t); \  y(t)=T_b(t) ; \ z(t)= T_c(t + \phi), $ where
$a,b$ are coprime integers, $c$ is an integer, $T_n(t)$ is the Chebyshev polynomial
of degree $n$ and $\phi \in \R.$
Chebyshev knots are polynomial analogues of Lissajous knots, which admit
parametrizations of the form
$x=\cos (at ); \   y=\cos (bt + \phi) ; \  z=\cos (ct + \psi)$,
where $ a, b, c$ are pairwise
coprime integers.  These knots, first defined in \cite{BHJS}, have been studied by many
authors: V. F. R. Jones, J. Przytycki,
C. Lamm, J. Hoste and L. Zirbel \cite{JP,La,HZ}.
\begin{figure}[ht]
\begin{center}
\begin{tabular}{ccc}
{\scalebox{.3}{\includegraphics{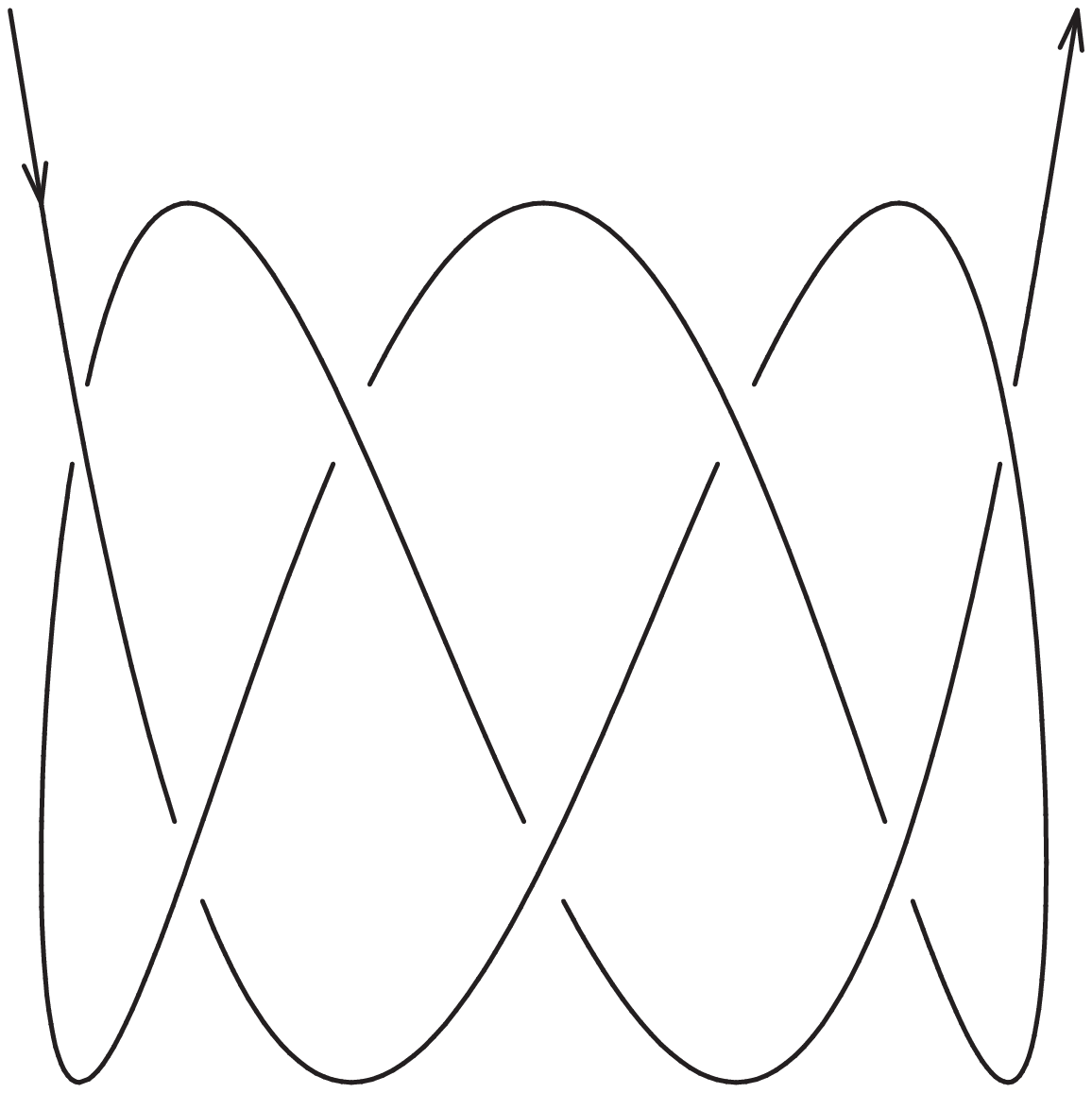}}}&\quad&
{\scalebox{.3}{\includegraphics{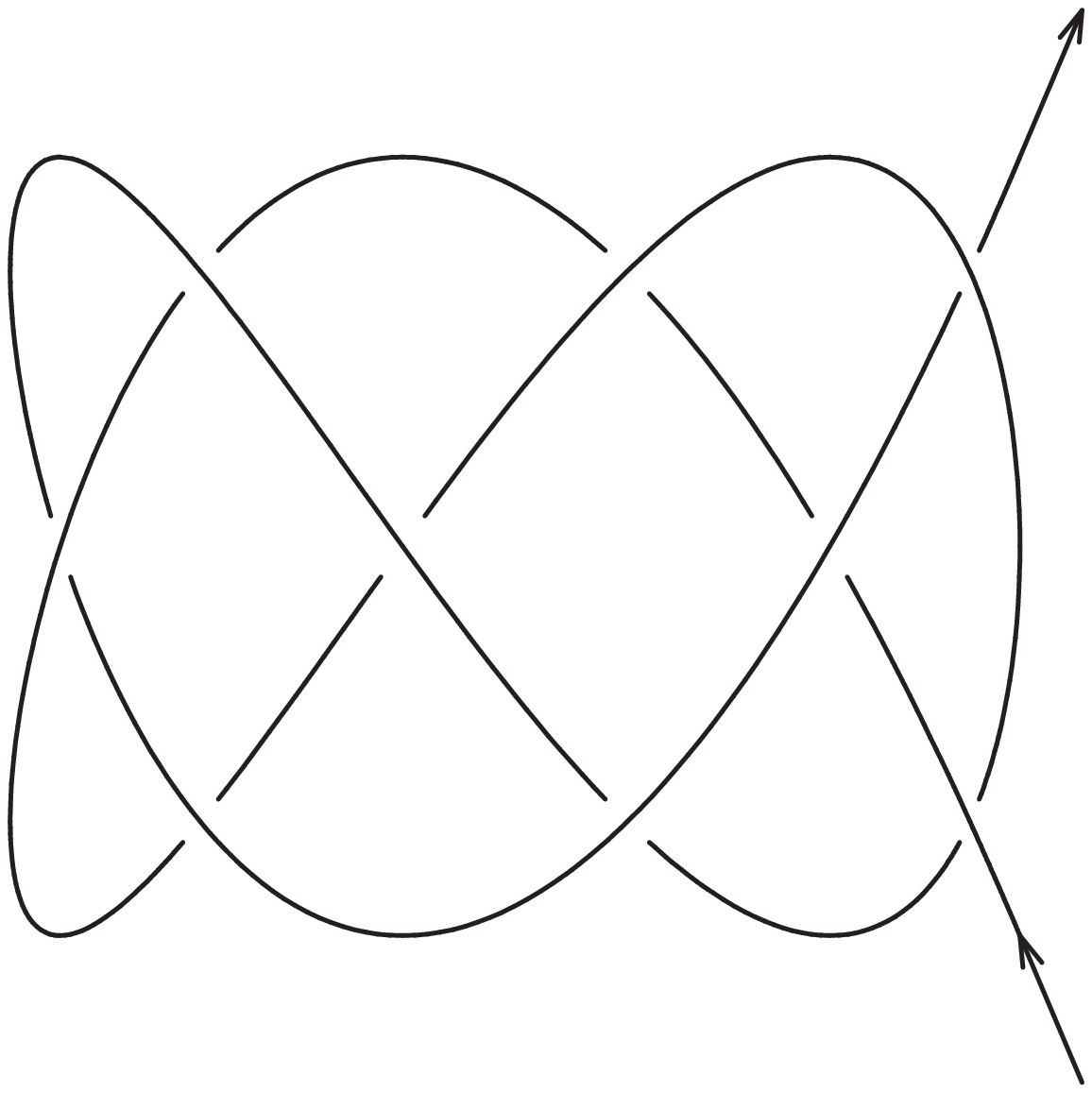}}}\\
$\cC(3,8,13,0)$ &\quad&
$\cC(4,7,9,0)$
\end{tabular}
\end{center}
\caption{Two Chebyshev knots: $\overline{7}_7$ and $7_5$ (Conway-Rolfsen numbering)}
\end{figure}
\pn
It is known that every knot may be obtained from a polynomial embedding
$\R \to \R^3 \subset \SS^3$, where $\SS^3$ is the one-point compactification of $\R^3$ (\cite{Va,DOS}).
In \cite{KP3} we have shown that every knot is a Chebyshev knot.

A two-bridge knot (or a rational knot) is a knot which is isotopic to a compact
space curve such that the $x$-coordinate has only two maxima and two minima.
When $a=3$ or $a=4$ the Chebyshev knot $\cC(a,b,c,\phi)$ is a two-bridge knot $K$.
Its projection onto
the $(x,y)$-plane is the Chebyshev curve $\cC(a,b): T_b(x)=T_a(y)$ and we have a regular diagram
of $K$ that is in Conway normal form when $a=3$ (see \cite{Con,Mu}).
This gives us an easy way to identify these knots using their classical Schubert invariant.
\pn
In \cite{KP4} we gave an explicit parametrization of the torus knots $T(2,n)$ and other
infinite families of knots. In \cite{KP4}, we gave a complete classification of harmonic
knots $\cC(a,b,c,0)$ where $a=3$ and $a=4$ (see also \cite{FF}). We have shown in \cite{KP4}
that every rational knot of crossing number $N$ admits a polynomial parametrization $x=T_3(t), \, y = T_b(t), \
z=C(t)$ with $b + \deg C = 3N$.
\pn
It is proved in \cite{KP3} that every two-bridge knot is a
Chebyshev knot with $a=3$. We first showed that every two-bridge knot has a Chebyshev
curve $\cC(3,b)$ as projection. Then we used a density argument based on Kronecker's theorem to
show that given a knot $K$ whose $\cC(a,b)$ is a plane projection, there exists $c$ and $\phi$ such that
$\cC(a,b,c,\phi)$ is isotopic to $K$.
\pn
Our aim is to give an algorithm that determines a Chebyshev parametrization $\cC(a,b,c,\phi)$,
(with $a=3$ and also with $a=4$), for any rational knot.
\bn
\item
Rational knots are classified by their Schubert fraction.
For $a=3$ and $a=4$ we determine the minimal integer $b$ such that the Chebyshev curve
$\cC(a,b): x=T_a(t),\, y=T_b(t)$ is a plane projection of $K$.
This algorithm is based on continued fraction expansions.
\item
Let $\cZ_{a,b,c}$ be the set of $\phi$ such that $\cC(a,b,c,\phi)$ is singular.
$\cZ_{a,b,c}$ is finite. The knot type of $K(\phi) = \cC(a,b,c,\phi)$ is constant over any
interval of $\R - \cZ_{a,b,c}$.
Then we determine a rational number in each component of $\R - \cZ_{a,b,c}$.
\item We determine the Schubert fraction of the knot $\cC(a,b,c,r)$
by evaluating the (under/over) nature of the crossings.
This amounts to evaluating the signs of polynomials at the real solutions of a zero-dimensional system.
\en
\pn
These algorithms are based on three basic black boxes:
\pn
\begin{itemize}
\item {\bf PhiProjection}($P_a(S,T), P_b(S,T), Q_c(S,T,\phi) \in \Q[\phi]$) \\
$\rightarrow$ $R \in \Q[\phi]$, such that $Z(R) = \pi_{\phi} \Bigl(\{ P_a=0,\ P_b=0, \ Q_c=0\}\Bigr)$
\item {\bf PhiSampling}($R \in \Q[\phi]$) \\
$\rightarrow$ $r_0, \ldots, r_N \in \Q$ such that $r_0<\phi_1<r_1<\cdots<\phi_N<r_N$, where
$\phi_1<\cdots<\phi_N$ are the real roots of $R$
\item {\bf SignSolve}($r \in \Q, \ P_a(S,T), \ P_b(S,T), D_{a,b,c}(S,T,\phi) \in Q[S, T,\phi]$) \\
$\rightarrow$ $[\sign {D_{a,b,c}(S,T,r)}], (S,T) \in \{P_a(S,T)=P_b(S,T)=0\}$
\end{itemize}
\pn
As a conclusion we give a list of Chebyshev parametrizations for the first 95 rational knots
of crossing number not greater than 10.
They admit polynomial parametrizations whose plane projections have
few crossing points. This is not the case with Lissajous knots (see \cite{BDHZ}).
\section{Geometry of Chebyshev knots}
Chebyshev curves were defined in \cite{Fi}.
Their double points are easier to study than those of Lissajous
curves. The classical Chebyshev polynomials $T_n$  are defined by
$T_n(t)= \cos  n \theta$, where $\cos \theta = t$.
These polynomials satisfy the linear recurrence
$T_0=1, \, T_1=t, \, T_{n+1} = 2t \, T_n - T_{n-1},\, n \geq 1$,
from which we deduce that
$T_n$ is a polynomial of degree $n$ and leading coefficient $2^{n-1}$.
\begin{proposition}[\cite{KP3}]\label{cab}
Let $a$ and $b$ be relatively prime integers.
The affine Chebyshev curve ${\cC(a,b)}$ defined  by
$$ {\cC(a,b)} \, : \quad T_b(x) - T_a(y) = 0$$
admits the parametrization $ x= T_a(t), \  y=T_b(t).$
$\cC$ has $\frac 12 (a-1)(b-1)$ singular points which are crossing points.
The pairs $(t,s)$ giving a crossing point are
$$t=\cos \left ( \Frac ka + \Frac hb
\right ) \pi,  \
s=\cos \left ( \Frac ka - \Frac hb \right ) \pi, $$
where $k,h$ are positive integers such that $\Frac ka + \Frac hb <1.$
\end{proposition}
When $a$ and $b$ are coprime integers, the projection
of $\cC(a,b,c,\phi)$ onto the $(x,y)$-plane is the plane Chebyshev curve
$\cC(a,b)$. The curve $\cC(a,b,c,\phi)$ is a knot if and only if $\cC(a,b,c,\phi)$ has no double point.
\pn
We thus deduce
\begin{proposition}\label{ck}
Let $a$, $b$ and $c$ be integers ($a$ and $b$ being relatively
prime). The number of  Chebyshev knots $\cC(a,b,c,\phi)$ is at most
$\frac 12 (a-1) \times (b-1) \times (c-1) + 1$.
\end{proposition}
\Pf
$\cC(a,b,c,\phi)$ is a singular space curve iff there exists
$s \not = t \in \R$ such that
$$
\Frac{T_a(t)-T_a(s)}{t-s} =0,\,
\Frac{T_b(t)-T_b(s)}{t-s} =0,\,
\Frac{T_c(t+\phi)-T_c(s+\phi)}{t-s} =0.
$$
From Prop. \ref{cab}, the set $\left\{\{s,t\}, \
\Frac{T_a(t)-T_a(s)}{t-s} =0,\,
\Frac{T_b(t)-T_b(s)}{t-s} =0 \right \}$
has $\frac 12(a-1)(b-1)$ elements.
For each of these elements, the set
$\{ \phi \in \R, \, \Frac{T_c(t+\phi)-T_c(s+\phi)}{t-s} = 0 \}$ has
at most $c-1$ elements because
the leading monomial of $\left [\Frac{T_c(t+\phi)-T_c(s+\phi)}{t-s} \right ]$
is $c2^{c-1} \phi^{c-1}$.
Consequently the set $\cZ_{a,b,c}$ of {\em critical values} $\phi$ has at most
$\frac 12 (a-1) \times (b-1) \times (c-1)$ elements.
For any $\phi \in \R - \cZ_{a,b,c}$, the curve $\cC(a,b,c,\phi)$ defines a knot.
We claim that in any interval included in $\R - \cZ_{a,b,c}$, the knots are
the same because the nature of the crossings is constant.
\EPf\pn
\begin{remark}
We see that $\cC(a,b,c,-\phi)$ is the reverse of $\cC(a,b,c,\phi)$ (see \cite{Mu}).
If $a=3$ or $a=4$, they define the same knot. $\cZ_{a,b,c}$ is symmetrical about the
origine.
\end{remark}
\section{Knot diagrams}
We shall study the diagram of the curve $\cC(a,b,c,\phi)$, that is to say the plane projection
$\cC(a,b)$ onto the $(x,y)$-plane and the nature (under/over) of the crossings over the double points
of $\cC(a,b)$.
There are two cases of crossing: the right twist and the left twist (see \cite{Mu}, P. 178).
\begin{figure}[th]
\begin{center}
\begin{tabular}{ccc}
{\scalebox{.2}{\includegraphics{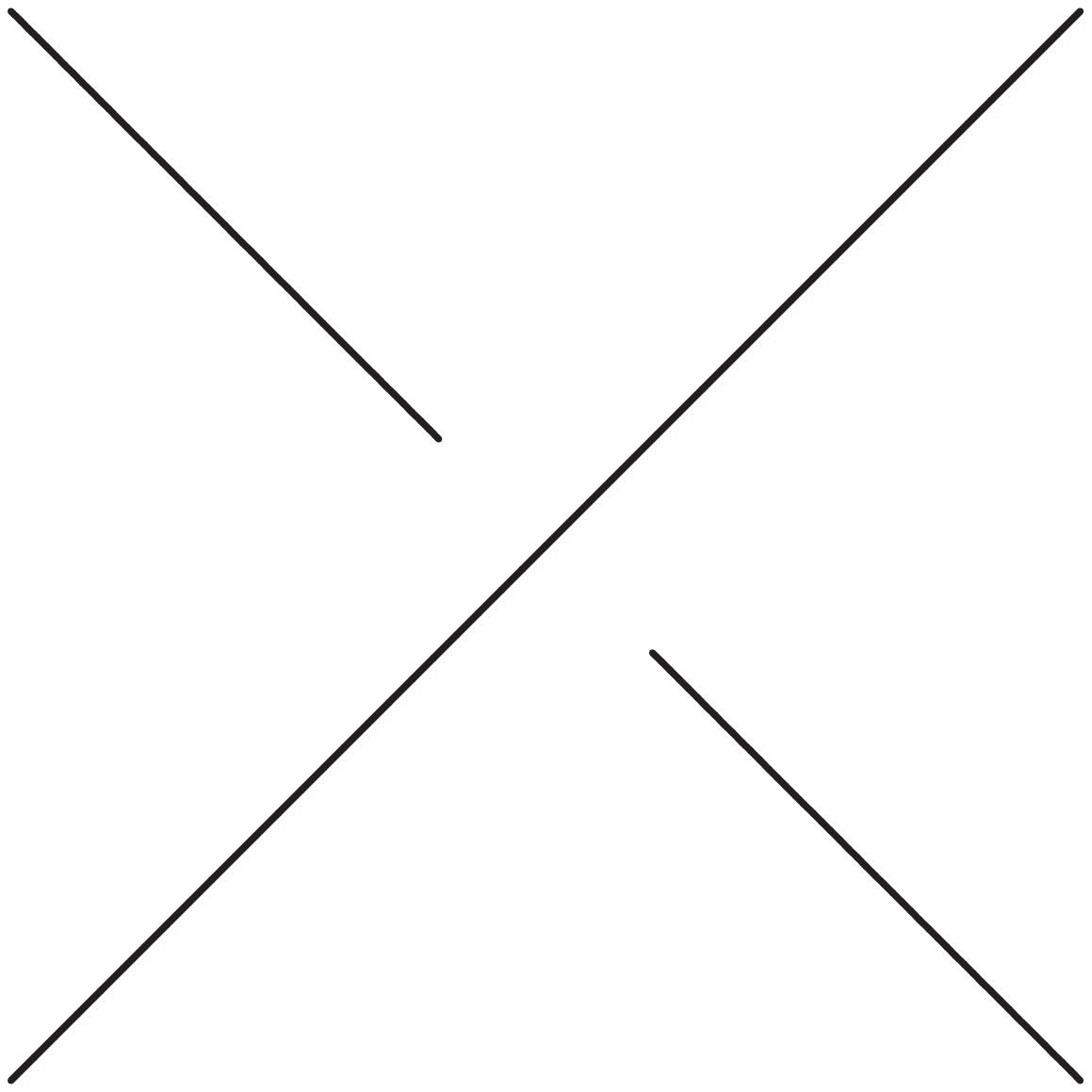}}} &\quad&
{\scalebox{.2}{\includegraphics{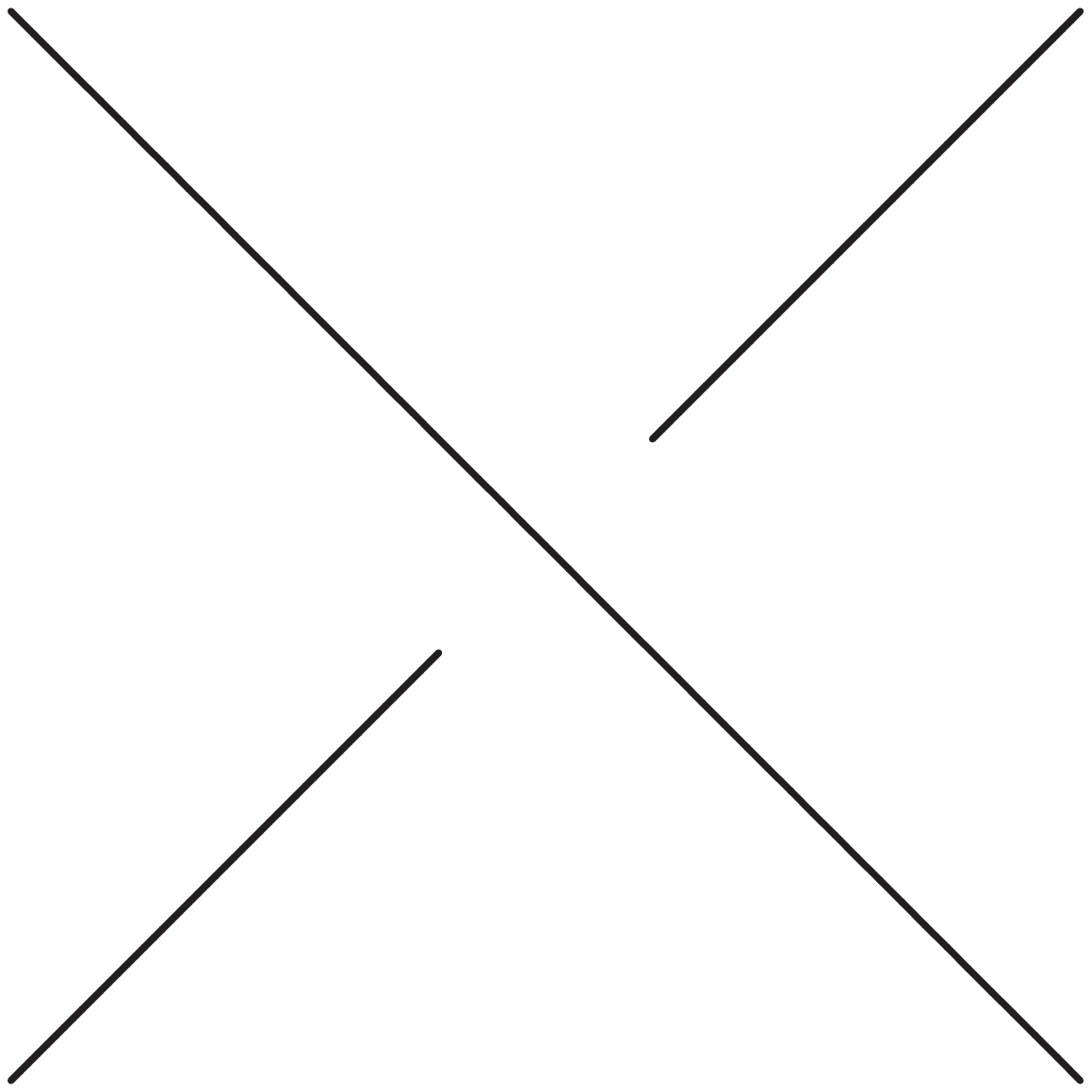}}}
\end{tabular}
\caption{The right twist and the left twist\label{signf}}
\end{center}
\end{figure}
In the following lemma, we see that the nature of the crossing is given by the sign of
a symmetrical polynomial.
\begin{lemma}
Let $s \not = t$ be parameters such that $T_a(t)=T_a(s)$ and $T_b(t)=T_b(s)$ and consider
the diagram of the curve $\cC(a,b,c,\phi)$. Let
\[
D(s,t,\phi) = \Bigl( \Frac{T_c(t+\phi)-T_c(t+\phi)}{t-s} \Bigr)
\Bigl(\Frac{T_{b-a}(t)-T_{b-a}(s)}{t-s}\Bigr). \label{D}
\]
Then $D(s,t,\phi)>0$ if and only if the crossing is a right twist.
\end{lemma}
\Pf
Let $(s,t)$ be the parameters of a double point of $\cC(a,b)$.
The crossing is a right twist if and only if
$$
\Bigl( z(t)-z(s)\Bigr)
\Bigl(\Frac{y'(t)}{x'(t)} - \Frac{y'(s)}{x'(s)}\Bigr)>0.
$$
Using Prop. \ref{cab}, we get $s=\cos \sigma$ and $t=\cos \tau$.
By simple computation we get
$\sin \tau \ x'(t) = - \sin \sigma \ x'(s)$, $\sin \tau \ y'(t) = \sin \sigma \ y'(s)$ and
$\Frac{y'(t)}{x'(t)} = - \Frac{y'(s)}{x'(s)}$. The slopes of the corresponding tangents
are opposite.
We therefore deduce that (using  $x\sim y$ for $\sign x = \sign y$)
$${x'(t)y'(t)} \sim - x'(s)y'(s) \sim \Frac{y'(t)}{x'(t)} - \Frac{y'(s)}{x'(s)}
\sim x'(t)y'(t) - x'(s)y'(s).$$
But
$$ x'(t)y'(t) = ab \Frac{\sin a \tau \ \sin b \tau}{\sin^2 \tau}
\sim 2 \sin a \tau \ \sin b \tau
= T_{a-b}(t) - T_{a+b}(t).$$
Consequently $x'(t)y'(t) - x'(s)y'(s) \sim (T_{a-b}(t) - T_{a+b}(t) ) - (T_{b-a}(s)-T_{b-a}(s)).$
On the other hand, using the identities $T_{b+a}+T_{b-a} = 2 T_aT_b$, $T_a(t)=T_a(s)$ and $T_b(t)=T_b(s)$,
we conclude that $x'(t)y'(t) - x'(s)y'(s) \sim T_{b-a}(t)-T_{b-a}(s)$ and the announced result.
\EPf\pn
A two-bridge knot (or link) admits a diagram in Conway's normal form.
This form, denoted by
$C(a_1, a_2, \ldots ,a_n)$  where $a_i$ are integers, is explained by the following
picture (see \cite{Con}, \cite{Mu} p. 187).
\psfrag{a}{\small $a_1$}\psfrag{b}{\small $a_2$}%
\psfrag{c}{\small $a_{n-1}$}\psfrag{d}{\small $a_{n}$}%
\begin{figure}[th]
\begin{center}
{\scalebox{.8}{\includegraphics{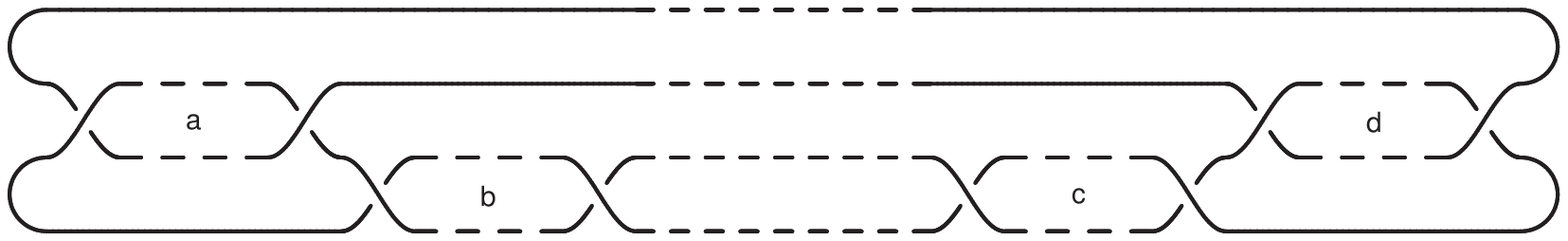}}}\\[30pt]
{\scalebox{.8}{\includegraphics{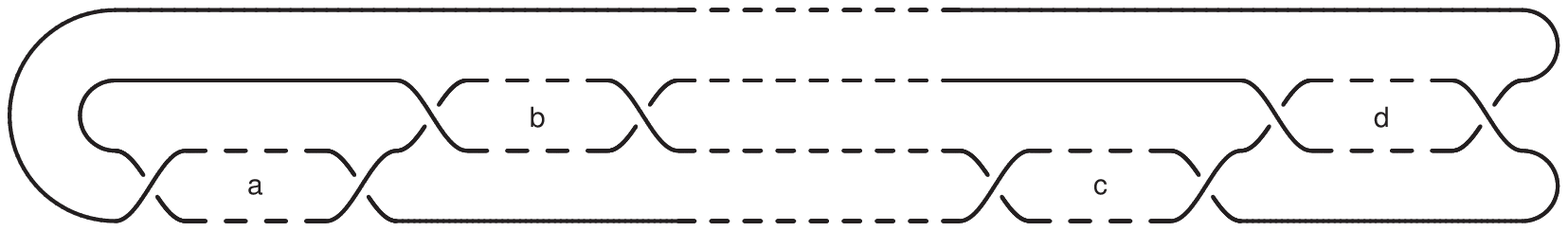}}}
\end{center}
\caption{Conway's normal forms, $n$ odd, $n$ even}
\label{conways3}
\end{figure}
The number of twists is denoted by the integer
$\abs{a_i}$, and the sign of $a_i$ is defined
as follows: if $i$ is odd, then the right twist is positive,
if $i$ is even, then the right twist is negative.
On  Figure \ref{conways3} the $a_i$ are positive (the $a_1$ first twists are right twists).
\pn
The two-bridge links are classified by their Schubert fractions
$$
\Frac {\alpha}{\beta} = [ a_1, \ldots , a_n], \quad \alpha >0,
$$
where $[ a_1, \ldots , a_n]$ is the continued fraction expansion
$a_1 + \Frac{1} {a_2 + \Frac {1} {a_3 + \Frac{1} {\cdots +\Frac 1{a_n}}}}$.

We shall denote by $S \bigl( \Frac {\alpha}{\beta} \bigr)$ the two-bridge link with
Schubert fraction $\Frac {\alpha}{\beta}.$
The two-bridge  links
$ S (\Frac {\alpha} {\beta} )$ and $ S( \Frac {\alpha ' }{\beta '} )$ are equivalent
if and only if $ \alpha = \alpha' $ and $ \beta' \equiv \beta ^{\pm 1} ( {\rm mod}  \  \alpha).$
The integer $ \alpha$ is odd for a knot, and even for a two-component link.
If $K= S (\Frac {\alpha}{\beta} ),$ its mirror image is
$ \overline{K}= S ( \Frac {\alpha}{- \beta} )$ (see \cite{Mu}).

We shall study knots with a Chebyshev diagram $\cC (3,b) : \  x= T_3(t), y= T_b(t).$
It is remarkable that such a diagram is already in Conway normal form
(see Figure \ref{chebd}).
Consequently, the Schubert fraction of such a knot is given by a
continued fraction of the form
$ [ \pm 1, \pm 1, \ldots ,\pm 1 ] .$
\pn
For example we obtain the torus knot $7_1 = C(-1,-1,-1,1,1,1,-1,-1,-1)$, and the knots
$6_3 = C(1,1,1,1,1,1)$, $6_1 = C(1,1,1,1,-1,-1,-1)$.
\def\p{{\small $+$}}
\def\m{{\small $-$}}
\begin{figure}[ht]
\begin{center}
\begin{tabular}{ccccc}
\psfrag{n1}{\m}\psfrag{n2}{\m}\psfrag{n3}{\m}\psfrag{n4}{\p}\psfrag{n5}{\p}\psfrag{n6}{\p}%
\psfrag{n7}{\m}\psfrag{n8}{\hspace{-3pt}\m}\psfrag{n9}{\hspace{-3pt}\m}%
{\scalebox{1.}{\includegraphics{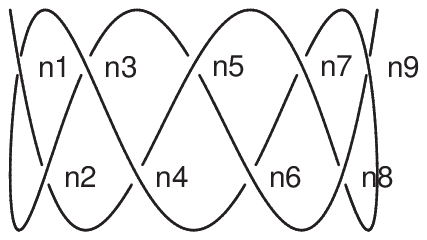}}}&\quad&
\psfrag{n1}{\p}\psfrag{n2}{\p}\psfrag{n3}{\p}\psfrag{n4}{\p}\psfrag{n5}{\p}\psfrag{n6}{\p}%
{\scalebox{1.}{\includegraphics{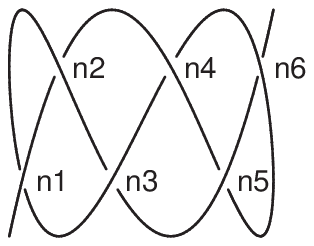}}}&\quad&
\psfrag{n1}{\m}\psfrag{n2}{\m}\psfrag{n3}{\m}\psfrag{n4}{\m}\psfrag{n5}{\p}\psfrag{n6}{\p}%
\psfrag{n7}{\p}%
{\scalebox{1.}{\includegraphics{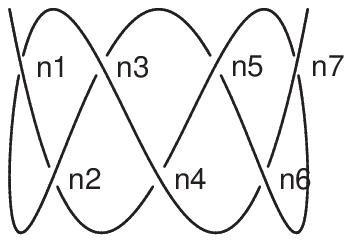}}}\\
$7_1$ & & $6_3$ && $\overline 6_1$
\end{tabular}
\caption{Chebyshev diagrams}\label{chebd}
\end{center}
\end{figure}
\par\noindent
We get for the knot $7_1$ (resp. $6_3$, $\overline{6_1}$) the fractions
$\Frac{7}{-6} \sim 7$, (resp. $\Frac{13}{8}$, $\Frac{9}{-5}$).

The {\em crossing number} of a knot is the smallest number of double points in
any plane projection of any isotopic knot.
The crossing number of a two-bridge knot $S(\Frac\alpha\beta)$ is the
sum of the integers in the regular continued fraction expansion of
$\Frac\alpha\beta >1$.

\begin{figure}[ht]
\begin{center}
\begin{tabular}{ccc}
\psfrag{n1}{\m}\psfrag{n2}{\m}\psfrag{n3}{\m}\psfrag{n4}{\p}\psfrag{n5}{\p}\psfrag{n6}{\p}%
{\scalebox{1.}{\includegraphics{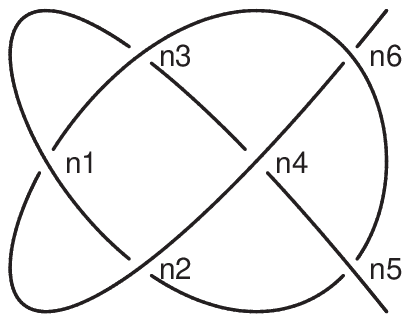}}}&\quad&
\psfrag{n1}{\p}\psfrag{n2}{\p}\psfrag{n3}{\p}\psfrag{n4}{\p}\psfrag{n5}{\p}\psfrag{n6}{\p}%
\psfrag{n7}{\p}\psfrag{n8}{\p}\psfrag{n9}{\p}%
{\scalebox{1.}{\includegraphics{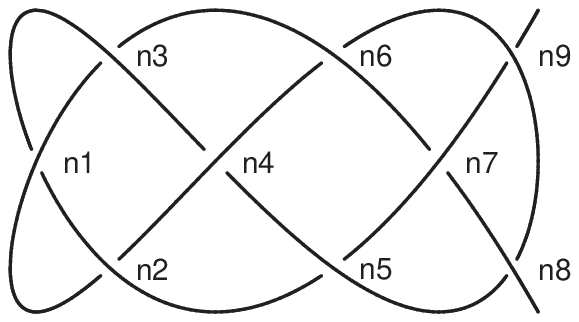}}}\\
$\overline{5}_2$ & & $9_{20}$
\end{tabular}
\caption{$a=4$}\label{a4}
\end{center}
\end{figure}

In the case of Chebyshev knots $\cC(3,b,c, \phi)$ with $a=3$, the
Conway notation is given by $C(\eps_1, \ldots, \eps_{b-1})$   where
$\eps_i = -(-1)^{i} \sign{D_i}$, the $D_i$ being ordered by their abscissae.
\pn
In the case of Chebyshev knots $\cC(4,b,c, \phi)$ with $a=4$, we obtain
diagrams like Figure \ref{a4}.
Let $a_i$ (resp. $b_i$, $c_i$)
be the signs of $D_i$ (resp. $-D_i$) corresponding to the crossing points
with $y=0$ (resp. $y<0$, $y>0$).
Following Murasugi (\cite{Mu}), the Conway normal form for such a knot is
$C(a_1, b_1+c_1, a_2, b_2+c_2, \ldots, a_n, b_n+c_n)$. Figure \ref{a4} shows the examples
$\overline{5}_2$: $C(-1,-2,1,2)$ and $9_{20}$: $C(1,2,1,2,1,2)$.
We thus deduce that the Conway notation for a knot $\cC(4,b,c,\phi)$ is
$C(a_1, b_1, \ldots, a_n, b_n)$ where $a_i = \pm 1$, $b_i = 0, \pm 2$.
\pn
In conclusion, we see that, in the particular case when $a=3$ or $a=4$,
the knot $\cC(a,b,c,r)$ is determined by its Schubert fraction $\Frac\alpha\beta$ corresponding
to the nature of the crossings over the double points of the projection $\cC(a,b)$.
\pn
On the other hand, we show that any rational number $\Frac\alpha\beta$ may be
expressed as continued fractions corresponding to Chebyshev diagrams $\cC(a,b)$ with $a=3$ and
$a=4$.
\begin{algorithm}\label{cfs}
Let $\Frac\alpha\beta$ be a rational number.
\bn
\item There exists a sequence
$\eps_1, \ldots, \eps_n$, $\eps_i=\pm 1$, such that
$\Frac\alpha\beta = [\eps_1, \ldots, \eps_n].$
\item If $\beta$ is even,
there exists a sequence
$\eps_1, \ldots, \eps_{2n}$, $\eps_{i}=\pm 1$,  such that\\
$\Frac\alpha\beta = [\eps_1, 2\eps_2, \ldots,\eps_{2n-1}, 2\eps_{2n}].$
\en
\end{algorithm}
\Pf
Let us prove the existence by induction on the height
$h(\Frac\alpha\beta)=\max(\abs{\alpha},\abs{\beta})$.
\bn
\item
\bi
\item[]
If $h=2$ then  $\Frac \alpha \beta = 1 = [1]$ and the  result is true.
\item[]
If $ \alpha > \beta,$ we have
$\Frac \alpha \beta = [1 ,\Frac \beta {\alpha - \beta }].
$
Since $h( \Frac \beta {\alpha - \beta } ) < h ( \Frac \alpha \beta )$,
we get our continued fraction by induction.
\item[]
If $ \beta > \alpha$  we have
$ \Frac \alpha \beta =
[1, -1, -\Frac {\beta - \alpha} \alpha  ].
$
And we also get the continued fraction.
\ei
This completes the construction of our continued fraction expansion $[\pm 1, \ldots, \pm1]$.
\item
\bi
\item[] If $h(\Frac\alpha\beta)=2,$
then $\alpha=1$ and $\beta=2$ and we have $r=[1,-2]$.
\item[] If $\alpha>2\beta>0$ then we write
$
\Frac\alpha\beta = [1,2,-1,2,\Frac{\alpha-2\beta}{\beta}]$.
We have $h(\Frac{\alpha-2\beta}{\beta}) < h(\Frac\alpha\beta)$
and we conclude by induction.
\item[] If $\beta<\alpha<2\beta$ then we write
$
\Frac\alpha\beta = [1,2,\Frac{\alpha-\beta}{3\beta-2\alpha}]$.
We have $\abs{3\beta-2\alpha}\leq \alpha$ and $\abs{\alpha-\beta}<\alpha$ and
we conclude by induction.
\item[] If $\beta > \alpha >0$ we write
$\Frac\alpha\beta = [1,-2,\Frac{\alpha-\beta}{2\alpha-\beta}].$
From $\abs{2\alpha-\beta}\leq \beta$ we have
$h(\Frac{\alpha-\beta}{2\alpha-\beta})<h(\Frac\alpha\beta)$ and we conclude
by induction.
\ei
\en
\hspace{1.truecm}
The existence of a continued fraction $[ 1, \pm 2, \ldots, \pm 1,  \pm2]$ is proved.
\EPf\pn
Note that we have proved in \cite{KP4} that the continued fraction expansion
$\Frac\alpha\beta = [\eps_1, \ldots, \eps_n]$, $\eps_i =\pm 1$, is unique if there
is no two consecutive sign changes and $\eps_{n-1}\eps_n>0$. We also proved that
the continued fraction expansion
$\Frac\alpha\beta = [\eps_1, 2\eps_2,\ldots,\eps_{2n-1}, 2\eps_{2n}]$, $\eps_i =\pm 1$,
is unique if there is no three consecutive sign changes.
\begin{cor}
Every two-bridge knot has a Chebyshev diagram $\cC(3,b)$, $b\not\equiv 0 \Mod 3$.
Every two-bridge knot has a Chebyshev diagram $\cC(4,b)$, $b \equiv 1 \Mod 2$.
\end{cor}
\Pf
Let us consider a knot $K= S(\Frac\alpha\beta)$.
Using Algorithm \ref{cfs}, we can write
$\Frac\alpha\beta = [\eps_1, \ldots, \eps_n], \eps_i = \pm 1$.
One can see that (see \cite{KP4}) $n \equiv 2\Mod 3$ iff $\alpha$ is even and, since $K$ is a knot,
this is not the case.
$K$ is isotopic to $C(\eps_1, \ldots, \eps_n)$ which
corresponds to a Chebyshev diagram $\cC(3,n+1): x=T_3(t), y=T_{n+1}(t)$.

Using Algorithm \ref{cfs}, we can write
$\Frac\alpha\beta = [\eps_1, 2\eps_2, \ldots,\eps_{2n-1}, 2\eps_{2n}], \eps_i = \pm 1$.
The knot $K$ is isotopic to $C(\eps_1, 2\eps_2, \ldots,\eps_{2n-1}, 2\eps_{2n})$ which
corresponds to a Chebyshev diagram $\cC(4,2n+1): x=T_3(t), y=T_{2n+1}(t)$.
\EPf\pn
\begin{cor}
Every two-bridge knot is a Chebyshev knot $\cC(3,b,c,\phi)$.
Every two-bridge knot is a Chebyshev knot $\cC(4,b,c,\phi)$.
\end{cor}
\Pf
Using a density argument (Kronecker theorem), we proved in \cite{KP3} that if $\phi$ is small enough,
then there exists $c$ such that $\cC(3,b,c,\phi) = C(\eps_1, \ldots, \eps_n)$.
The case $a=4$ is similar.
\EPf\pn
Unfortunately, this last corollary will not provide $c$ and $\phi$ and not even any bound for $c$.
We want to give the minimal Chebyshev parametrization for every rational knot with a small crossing
number. We shall describe all rational knots $\cC(a,b,c,\phi)$ with given $a, b, c$.
\section{Description of Chebyshev knots}\label{descr}
Let us consider the curve $\cC(a,b,c,\phi)$ with $a=3$ or $a=4$. From section 3., we know that
\bn
\item The curve is singular iff it has double points.
\item If the curve is not singular, the knot $\cC(a,b,c,\phi)$ is determined by the sequence
of crossings of the projection $\cC(a,b)$.
\en
\pn
We will use the symmetric variables $S=s+t$ and $T=st$. Let us define
\[P_n (S,T) = \Frac{T_n(t)-T_n(s)}{(t-s)}, \
Q_n(S,T,\phi) = \Frac{T_n(t+\phi)-T_n(s+\phi)}{(t-s)}.\label{PnQn}
\]
\begin{lemma}
There exists $R_{a,b,c} \in \Q[\phi]$ with degree
$\deg R_{a,b,c} \leq \frac 12 (a-1)(b-1)(c-1)$ such that $\cZ_{a,b,c}=Z(R_{a,b,c})$.
\end{lemma}
\Pf
From Prop. \ref{cab}, $\{ P_a(S,T)=0, \, P_b(S,T)=0 \}$ is 0-dimensional and has degree
$\frac 12 (a-1)(b-1)$.
$V_{a,b,c} = \{ P_a(S,T)=0, \, P_b(S,T)=0, \, Q_c(S,T,\phi)=0\}$
is 0-dimensional and has degree $\frac 12 (a-1)(b-1)(c-1)$ because
$Q_c$ has $c (2\phi)^{c-1}$ as leading term.
From lemma \ref{ck}, $\cZ_{a,b,c} = \pi_{\phi}(V_{a,b,c})$ is finite. It is for example
$Z(R)$ where $\langle R \rangle = V_{a,b,c} \bigcap \Q[\phi]$.
\EPf\pn
We obtain $R_{a,b,c}$ by using the black-box ${\tt PhiProjection}(P_a,P_b,Q_c)$.
\pn
The set $\cZ_{a,b,c} = \{\phi_1, \ldots, \phi_N\}$,
for which $\cC(a,b,c,\phi)$ is singular, is exactly
$\pi_{\phi} (V_{a,b,c})$, from Prop. \ref{ck}. 
From section 3., the knot $K(r) = \cC(a,b,c,r)$ is constant over any interval
$]\phi_i, \phi_{i+1}[$. 
We obtain $r_i$ in $]\phi_i, \phi_{i+1}[$ by using ${\tt PhiSampling}(R_{a,b,c})$.
\pn
In the case when $a=3$ or $a=4$, and $r \in \R - \cZ_{a,b,c}$, the knot
$K(r)$ is uniquely determined by its Schubert fraction.
Let $(S_i,T_i)_{i=1,\ldots,n}$ be the points of $\{P_a(S,T)=0,\, P_b(S,T)=0\}$.
The Conway notation of the knot $K(r)$ is deduced from
the sequence $[\sign{D(S_i,T_i,r)}]_{i=1,\ldots,n} = {\tt SignSolve}(r,P_a, P_b, Q_c\cdot P_{b-a})$.
\pn
We will now show how to determine these quantities in the case $a=3$ and $a=4$.
\subsection{Case ${a=3}$}\label{c3}
We get $P_3(S,T) = 4T-4S^2+3$ that is $T=S^2 - \frac 34$.
The set $\{(S,T),\, P_3(S,T)=0,P_b(S,T)=0\}$ has cardinal $b-1$. Its elements satisfy
$P_b(S,S^2 - \frac 34)=0, \, T=S^2 - \frac 34$. We deduce that $\deg(P_b(S,S^2 - \frac 34))=b-1$.
The set of critical values $\phi$ is
$$
\cZ_{3,b,c}=\pi_{\phi} \Bigl(\{P_b(S,S^2 - \frac 34)=0, \, Q_c(S,S^2 - \frac 34,\phi)=0\}\Bigr).
$$
It is exactly the roots of the polynomial of degree $(b-1)(c-1)$:
$$R_{3,b,c} = \Res_S\Bigl(P_b(S,S^2 - \frac 34),Q_c(S,S^2 - \frac 34,\phi)\Bigr).$$
Let $A(S)$ be a crossing point corresponding to parameter $S$ (and $T=S^2-\frac 34$).
Its abscissa is $T_3(t)=T_3(s)= -T_3(S) = S(3-4S^2)$.
\psfrag{a1}{}\psfrag{b0}{\small $A_{b-2}$}\psfrag{c0}{}%
\psfrag{a0}{\small $A_{b-1}$}%
\psfrag{cn}{\small $A_{1}$}\psfrag{bn}{\small $A_{2}$}\psfrag{an}{\small $A_{3}$}%
\begin{figure}[th]
\begin{center}
{\scalebox{.7}{\includegraphics{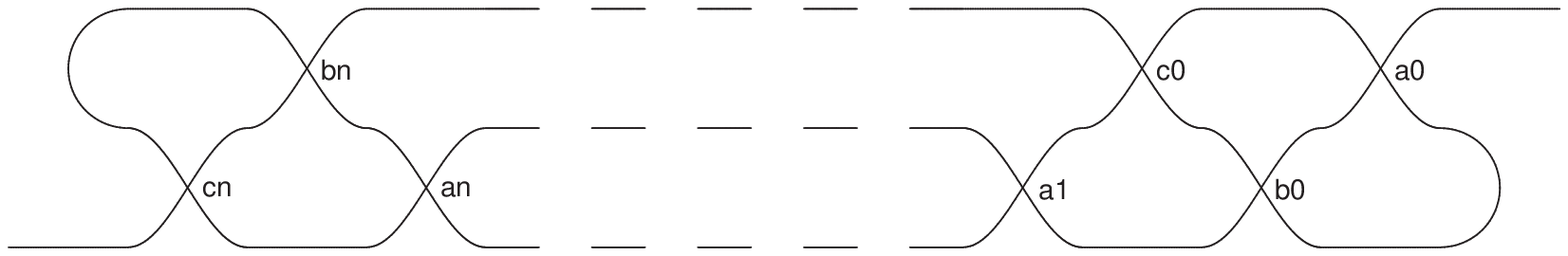}}}
\end{center}
\caption{$\cC(3,b)$, $b$ even}
\label{dh3}
\end{figure}
\pn
We define the order relation $A(S)<_3 A(S')$ if $T_3(S)>T_3(S')$.
The Conway notation of $\cC(3,b,c,r)$ is
$C(D(A_1),-D(A_2),\ldots, (-1)^{b-1}D(A_{b-1}))$ where $A_1 <_3 A_2 <_3 \cdots <_3 A_{b-1}$ and $D(A)$ is $D(s,t,r)$ defined in Formula \ref{D}.
\subsection{Case ${a=4}$}\label{c4}
We get $P_4(S,T) = 8\,S \left({S}^{2} -2\,T-1 \right)$.
We thus obtain two families of double points
$$\cA: \{ S=0,\, P_b(0,T) = 0\}, \quad
\cB: \{T=\Frac 12 \left ( S^2-1 \right ), P_b(S,\frac 12 ( S^2-1 )) = 0\}.
$$
Let $n =\frac 12 (b-1)$.
We have $\abs{{\cal A}} = n$ from which we deduce that
$\deg_T P_b(0,T) = n$. From $\abs{{\cal B}} = 2n$, we deduce also
that $\deg_S P_b(S,\frac 12 ( S^2-1 )) = 2n$.
As the leading coefficient of $Q_c(S,T,\phi)$ is $c(2\phi)^{c-1}$ we deduce
that
$R_1(\phi) = \Res_T(P_b(0,T), Q_c(0,T,\phi))$
has degree $n(c-1)$ and
$R_2(\phi) = \Res_S\Bigl(P_b(S,\Frac 12 \left ( S^2-1 \right )), Q_c(S,\Frac 12 \left ( S^2-1 \right ),\phi)\Bigr)$
has degree $2n(c-1)$.
$\cZ_{4,b,c}$ is the set of real roots of $R_{4,b,c}= R_1 \times R_2$.
\psfrag{a0}{$B_{n}$}%
\psfrag{aa0}{$B'_{n}$}%
\psfrag{c0}{$B_{n-1}$}%
\psfrag{cc0}{$B'_{n-1}$}%
\psfrag{d0}{$A_{n-1}$}%
\psfrag{b0}{$A_{n}$}%
\psfrag{cn}{$B_{1}$}%
\psfrag{dn}{$A_{1}$}%
\psfrag{ccn}{$B'_{1}$}%
\psfrag{bn}{$A_{2}$}%
\psfrag{an}{$B_{2}$}%
\psfrag{aan}{$B'_{2}$}%
\begin{figure}[th]
\begin{center}
{\scalebox{.7}{\includegraphics{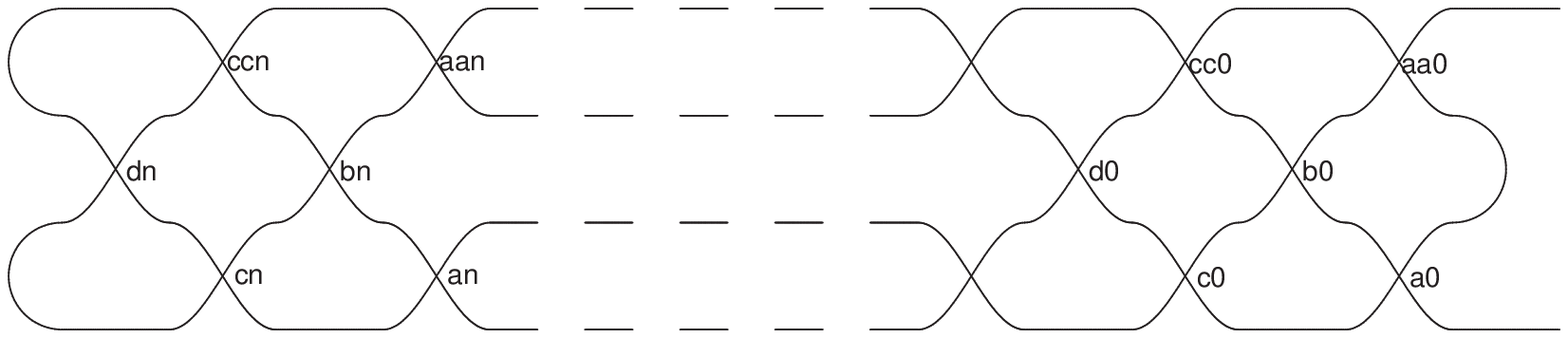}}}
\end{center}
\caption{$\cC(4,2n+1)$}
\label{dh4}
\end{figure}
\pn
The abscissa of $A(T) \in \cA$ is given by $T_4(t)=T_4(s)= 1+8\,T+8\,T^2.$
The abscissa of $B(S) \in \cB$ is given by $T_4(t)=T_4(s)=-1+4\,S^2-2\,S^4$.
We have to sort separately the crossing points of $\cA$ and $\cB$ by increasing abscissae:
$A_1, \ldots, A_{n}$ and $B_1,B'_1, \ldots, B_{n},B'_{n}$. Note
that $B_i$ and $B'_i$ have the same abscissa.
The Conway notation for the knot we obtain is then
$$C\Bigl(D(A_1),-( D(B_1)+D(B'_1)), \ldots, D(A_{n}), -( D(B_{n})+D(B'_{n}))\Bigr),$$
where $D(A)$ (resp. $D(B)$) is $D(s,t,r)$ defined in Formula \ref{D}.
\begin{remark}
We have here $\deg R_{a,b,c} = \frac 12(a-1)(b-1)(c-1)$. We could have computed also
$R_{a,b,c}$ by eliminating $S$ and $T$ using Gr\"obner Basis (see \cite{CLOS}).
It may happen that $\deg R_{a,b,c} < \frac 12(a-1)(b-1)(c-1)$.
\end{remark}
\subsection{Computation of the polynomials}
As $T_n$ satisfies the linear recurrence of order 2: $T_{n+1}+T_{n-1} = 2t T_n$ we deduce that
$Q_n$ (Form. \ref{PnQn}) satisfies the linear recurrence of order 4:
\[
&Q_0 = 0, \, Q_1 = 1, \,
Q_2 = 2S+ 4 \phi, \, Q_3 = -4\,T+12\,\phi\,S+4\,{S}^{2}+12\,{\phi}^{2}-3.&\nonumber\\
&Q_{n+4} = 2 \left(S + 2\,\phi \right)\left (Q_{n+3}+Q_{n+1} \right )
- 2 \left( 2\,{\phi}^{2}+ 2\,T+2\,\phi\,S + 1\right) Q_{n+2} - Q_c.& \label{Q_n}
\]
For $P_n(S,T) = Q_n(S,T,0)$ we find
\[
&P_0 = 0, \, P_1 = 1, \,
P_2 = 2S, \, P_3 = -4\,T+4\,{S}^{2}-3.&\nonumber\\
&P_{n+4} = 2 S\left (P_{n+3}+P_{n+1} \right )
- \left(4\,T+2 \right) P_{n+2} - P_n.& \label{P_b}
\]
In the particular case when $a=3$ or $a=4$ we have to compute
$Q_n(S,T,\phi)$ where $T=S^2 -\frac 34$ or $T=\frac 12 (S^2-1)$ or $S=0$. These polynomials satisfy
also linear recurrences.

As $T_m(T_n) = T_{mn}$ we deduce that $Q_n \vert Q_{nm}$. We can therefore obtain factors
of the polynomials $Q_n(S,T,\phi)$.

In the particular case where $a=3$ or $a=4$, we obtain our resultants
$R_{a,b,c} = {\tt PhiProjection}(P_a,P_b,Q_c)$ by computing the resultants between factors
of $P_b(S,T)$ and factors of $Q_c(S,T,\phi)$ that depend only on $S$ and $\phi$ or
on $T$ and $\phi$.
\pn
We have to determine the Schubert fraction of any knot of the
type $\cC(a,b,c,r)$ where $r$ is a given rational number in $\Q - \cZ_{a,b,c}$. Such
rational number is given by ${\tt PhiSampling}(R_{a,b,c})$.
\section{Examples}
\subsection{The family of knots $\cC(3,5,7,\phi)$}
We get $P_3(S,T) = \Frac{T_3(t)-T_3(s)}{t-s} = 4(T-S^2+\frac 34)$, so $T=S^2-\frac 34$ and
$$
\begin{array}{rcl}
P_5(S,S^2-\frac 34) &=& -16 S^4+12 S^2-1,\\
Q_7(S,S^2-\frac 34,\phi) &=&
64\,{S}^{6}-16\, \left( 84\,{\phi}^{2}+5 \right) {S}^{4}-112\,\phi\,
\left( 1+20\,{\phi}^{2} \right) {S}^{3}+24\, \left( 42\,{\phi}^{2}+1
\right) {S}^{2}\\
&&\quad +28\,\phi\, \left( 48\,{\phi}^{4}+3+80\,{\phi}^{2}
\right) S-1+1120\,{\phi}^{4}+448\,{\phi}^{6}+84\,{\phi}^{2}.
\end{array}
$$
\psfrag{S1}{\small $S_1$}\psfrag{S2}{\small $S_2$}%
\psfrag{S3}{\small $S_3$}\psfrag{S4}{\small $S_4$}%
\def\v1{\vspace{10pt}\par\noindent}
\def\h1{\hspace{-15pt}}
\psfrag{S1}{\small $S_1$}\psfrag{S2}{\small $S_2$}%
\psfrag{S3}{\small $S_3$}\psfrag{S4}{\small $S_4$}%
\psfrag{f1}{\small $\phi_1$}%
\psfrag{f2}{\small $\phi_2$}%
\psfrag{f3}{\small $\phi_3$}%
\psfrag{f4}{\h1\small $\phi_4$}%
\psfrag{f5}{\small $\phi_5$}%
\psfrag{f6}{\h1\small $\phi_6$}%
\psfrag{f7}{\small $\phi_7$}%
\psfrag{f8}{\h1\small $\phi_8$}%
\psfrag{f9}{\small $\phi_9$}%
\psfrag{f10}{\h1\small $\phi_{10}$}%
\psfrag{f11}{\small $\phi_{11}$}%
\psfrag{f12}{\small $\phi_{12}$}%
\begin{figure}[th]
\begin{center}
{\scalebox{.8}{\includegraphics{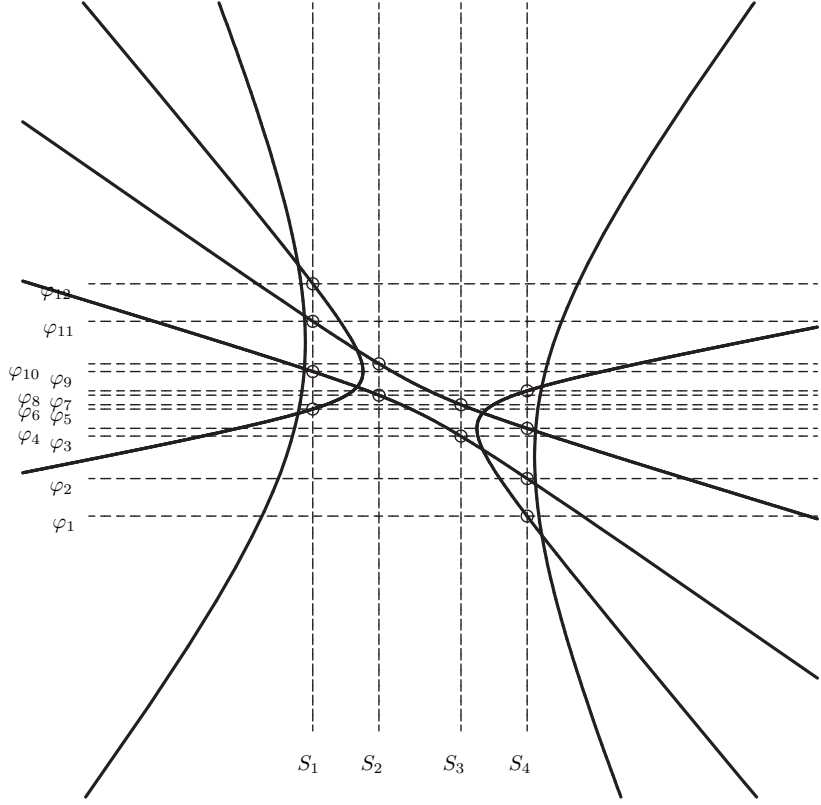}}}
\end{center}
\caption{The curve $Q_7(S,S^2-\frac 34,\phi) = 0$ and the lines $P_5(S,S^2-\frac 34) = 0$}
\label{r357}
\end{figure}
\pn
$\cC(3,5,7,\phi)$ is singular iff $\phi$ is a root of
$R_{3,5,7}=\Res_{S}\Bigl(P_5(S,S^2-\frac 34),Q_7(S,S^2-\frac 34, \phi)\Bigr)$.
\pn
$R_{3,5,7}$ has degree $24=\frac 12 (3-1)(5-1)(7-1)$ and 12 real roots $\phi_1, \ldots, \phi_{12}$.
We choose 13 rational values $r_0<\phi_1<r_1 < \cdots < \phi_{12}<r_{12}$.
\pn
Let us determine now the nature of $\cC(3,5,7,r)$.  We have to evaluate
$D(s,t,\phi) = Q_7(S,S^2-\frac 34,r) \cdot P_{2}(S,S^2-\frac 34)$
when $P_5(S,S^2-\frac 34)=0$. Let $S_1 < S_2< S_3 < S_4$ be the 4 real roots of $P_5$.
They correspond to parameters $(s_1,t_1), \ldots, (s_4,t_4)$ such that
$s_i+t_i = S_i, \, s_i t_i = T_i = S_i^2-\frac 34$.
We have $T_3(-S_2)<T_3(-S_1)<T_3(-S_4)<T_3(-S_3)$ and the knot $\cC(3,5,7,r)$ is given by the
continued fraction expansion $\alpha/\beta=[D_2,-D_1,D_4,-D_3]$ where $D_i = \sign {Q_7(S_i,S_i^2-\frac 34,r)\cdot
P_2(S_i,S_i^2-\frac 34)}$. We obtain
\setlength{\tabcolsep}{5pt}
\begin{center}
\begin{tabular}{||c|c|c|c|c|c|c||}
\hline
$r_7= 0$&$r_8= \frac 1{15}$&$r_9= \frac 15$&$r_{10}= \frac 14$&
$r_{11}= \frac 12$&$r_{12}= \frac 23$&$r_{13}= 1$\\
\hline
{\scalebox{1}{\includegraphics{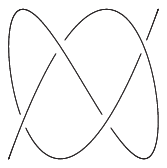}}}&
{\scalebox{1}{\includegraphics{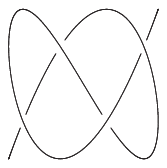}}}&
{\scalebox{1}{\includegraphics{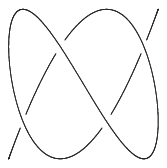}}}&
{\scalebox{1}{\includegraphics{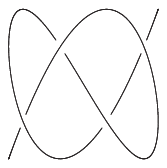}}}&
{\scalebox{1}{\includegraphics{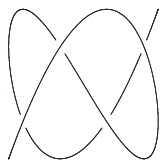}}}&
{\scalebox{1}{\includegraphics{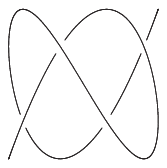}}}&
{\scalebox{1}{\includegraphics{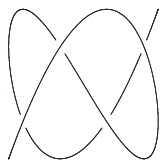}}}\rule[-7pt]{0pt}{2.2cm}\\
\hline
$\frac\alpha\beta = \frac 53$ &$\frac\alpha\beta =  -\frac {1}{3}$ &$\frac\alpha\beta =  -1$&$\frac\alpha\beta =  -1$ &$\frac\alpha\beta =  1 $
&$\frac\alpha\beta = 1$ &$\frac\alpha\beta = 1$\\
\hline
\end{tabular}
\end{center}
The only non trivial knot is obtained for $r=0$. It is the figure-eight knot $4_1 = S(\Frac 52)$.
\subsection{A more complicated example, the family of knots $\cC(3,14,292,\phi)$}\label{ex2}
$P_{14}(S,S^2-\frac 34)$ is a product of 4 factors of degrees $[6,3,3,1]$.
We have $292 = 4\cdot 73$. We know that $Q_{73} \vert Q_{146} \vert Q_{292}$.
$Q_{73}(S,S^2-\frac 34,\phi)$ is irreducible and has degree 72. $Q_{146}$ is a product of polynomials with
degrees $[72,72,1]$. At the end $Q_{292}$ is a product of 5 factors with degrees $[144,72,72,2,1]$.
We compute $R_{3,14,292}$ as the product of 20 resultants between factors of $P_{14}$ and
$Q_{292}$. $R_{3,14,292}(\phi)$ has degree 3783 and exactly 2185 distinct real roots. We
compute the 1093 Schubert fractions $\cC(3,14,292,r_i)$ where $r_i>0$.
We obtain 275 non trivial knots and eventually 34 distinct knots. One of these
has crossing number greater that 10, it is the knot
$\overline{12}_{518} = S(\frac{157}{34})$.
\subsection{A much more complicated example, the family of knots $\cC(4,13,267,\phi)$}
$P_{13}(S,T)$ is irreducible so as $P_{13}(0,T)$ (that has degree 6) and
$P_{13}(S,\frac 12(S^2-1))$ (that has degree 12).
We have $267 = 3\cdot 89$. We know that $Q_{89} \vert Q_{267}$. $Q_{267}$ is the product
of $Q_3$, $Q_{89}$ and a polynomial of degree 176 in $\phi$.
$Q_{267}(S,T,\phi)$ is a product of polynomials of degrees $[176,88,2]$.

We thus obtain $R_1 = \Res_T \Bigl(P_{13}(0,T),Q_{267}(0,T)\Bigr)$ as a product of
polynomials of degrees $[1056, 528, 12]$ and
$R_2 = \Res_S \Bigl(P_{13}(S,\frac 12(S^2-1)),Q_{267}(S,\frac 12(S^2-1))\Bigr)$ as
a product of polynomials of degrees $[2112, 1056, 24]$.

$R_{4,13,267} = R_1 \times R_2$ has degree
4788. It has 2882 distinct real roots. We compute the 1442 Schubert fractions
$\cC(4,13,267,r_i)$ where $r_i>0$.
We obtain 710 non trivial knots. 72 of these are distinct knots whose crossing numbers
take all values between 3 and 16.
\section{Results}
In this section we present some results we have obtained using certified implementations
of the three black-boxes on which our algorithms are based. They are
easily implementable in any high level language.
\subsection{Implementations}
There are numerous choices for the implementations, but our
requirements are strict: we must certify all the results since our goal is
to obtain a classification; bearing in mind that
the systems of polynomial equations have thousands of roots.
\pn
\begin{itemize}
\itemsep 3pt plus 1pt minus 1pt
\item {\bf PhiProjection}($P_a(S,T), P_b(S,T), Q_c(S,T,\phi) \in \Q[\phi]$) \\
$\rightarrow$ $R \in \Q[\phi]$, such that $Z(R) = \pi_{\phi} \Bigl(\{ P_a=0,\ P_b=0, \ Q_c=0\}\Bigr)$
\pn
A straightforward way is to compute a Gr\"obner basis of
$\langle P_a(S,T), \, P_b(S,T), \, Q_c(S,T,\phi) \rangle$ for a so called
{\em elimination order} (see \cite{CLOS}).
Triangular decompositions provide a suitable alternative, or,
more basically, iterative resultants in {\em generic} situations.
Resultants can be used efficiently for our
problem, the system being sufficiently generic. Our choice is then an {\em ad-hoc} method based
on resultants computations in the same spirit as in section \ref{descr}.
\item {\bf PhiSampling}($R \in \Q[\phi]$) \\
$\rightarrow$ $r_0, \ldots, r_N \in \Q$ such that $r_0<\phi_1<r_1<\cdots<\phi_N<r_N$, where
$\phi_1<\cdots<\phi_N$ are the real roots of $R$
\pn
Such a function can easily be implemented using any solver that is able to {\em isolate} real roots of
univariate polynomials (say providing non overlapping intervals with rational bounds around all the
real roots). It must be able to discriminate multiple roots from clusters of roots, real roots
from complex roots with a small imaginary part (which excludes many numerical methods and most of
implementations using hardware floats). One can use methods based on Sturm sequences or the Descartes
rule of signs (see \cite{BPR} for an overview), but also many strategies using interval analysis.
Due to the high degree of the polynomials, our choice is to use algorithms based on the Descartes rule
of signs using multi-precision interval arithmetic as in \cite{RZ}.
\item {\bf SignSolve}($r \in \Q, \ P_a(S,T), \ P_b(S,T), D_{a,b,c}(S,T,\phi) \in Q[S, T,\phi]$) \\
$\rightarrow$ $[\sign {D_{a,b,c}(S,T,r)}], (S,T) \in \{P_a(S,T)=P_b(S,T)=0\}$
\pn
The determination of the sign of a polynomial over a zero-dimensional system is
difficult to certify when using numerical method. There are few exact/certified existing
methods/implementations for this problem. The strategy is naturally linked
to the implementation of the function {\tt PhiProjection} since the
zero-dimensional system to be considered by {\tt SignSolve} is a subsystem
of the one which is to be considered by {\tt PhiProjection}.
One can use the {\em generalized Hermite method} for zero-dimensional systems
as in \cite{PRS}, which makes use of Gr\"obner bases.
One can use also any method that first rewrites the system as a
rational parametrization (as in \cite{Rou} or \cite{GLS}) and then apply any algorithm
that computes the sign of an univariate polynomial at a real algebraic number.
This last step can be done by extending methods based on Sturm theorem or based on the Descartes rule
of signs. Due to our implementation of {\tt PhiProjection} and to the degrees of the
polynomials, we base our implementation on the Descartes rule of signs.
\end{itemize}
For the experiments, we used the {\sc Maple} environment.
{\tt PhiProjection} is based on resultants computation (see \ref{c3} and \ref{c4}).
We use the {\sc Maple} function {\tt Isolate} for {\tt PhiSampling} (without constraints)
and for {\tt SignSolve} (with constraints). In the univariate case,
this function is based on the algorithm described in \cite{RZ}.
Other computations have been straightforwardly implemented
according to the descriptions proposed in section \ref{descr}.
\subsection{Experiments}
\begin{figure}[th]
\begin{center}
\begin{tabular}{ccc}
{\scalebox{.8}{\includegraphics{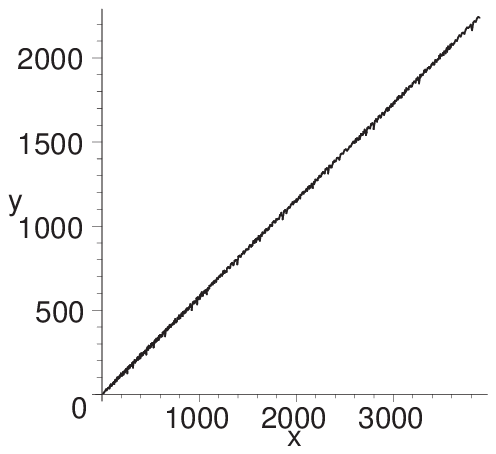}}}&
{\scalebox{.8}{\includegraphics{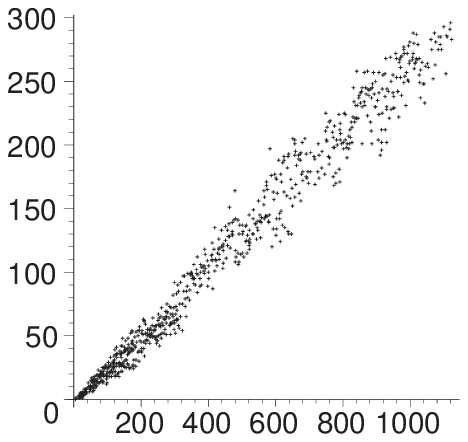}}}&
{\scalebox{.8}{\includegraphics{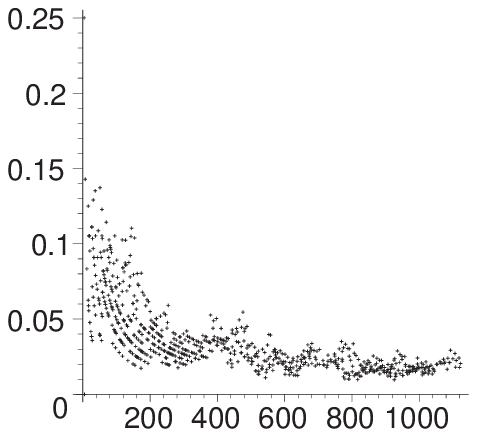}}}\\
$\abs{\cZ_{3,b,c}}$& Non trivial knots, $a=3$ & Distinct knots, $a=3$\\
{\scalebox{.8}{\includegraphics{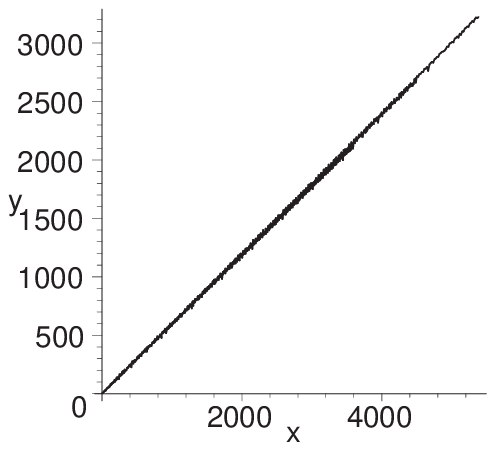}}}&
{\scalebox{.8}{\includegraphics{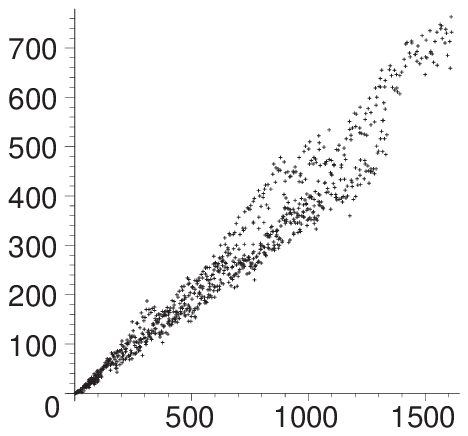}}}&
{\scalebox{.8}{\includegraphics{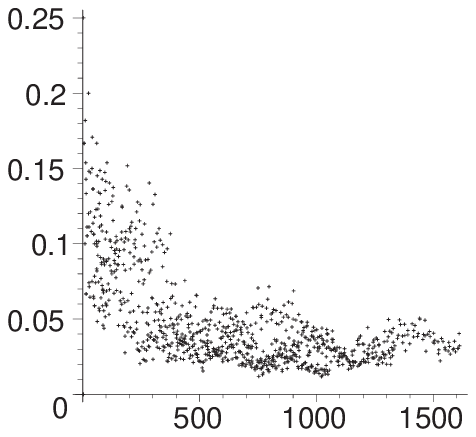}}}\\
$\abs{\cZ_{4,b,c}}$& Non trivial knots, $a=4$ & Distinct knots, $a=4$\\
\end{tabular}
\end{center}
\caption{Number of real roots, number of non trivial knots, proportion of distinct knots}
\label{racs}
\end{figure}

Let us remind in the next table the number $\cK_N$ of two-bridge knots with crossing number $N$,
up to mirror symmetry (see \cite{ES} for a formula).
\pn
\begin{center}
\begin{tabular}{||c|cccccccc||}
\hline
$N$&3&4&5&6&7&8&9&10\\
\hline
$\cK_N$  &1&1&2&3&7&12&24&45\\
\hline
\end{tabular}
\end{center}
\pn
The minimal $b$ for a Chebyshev diagram $\cC(3,b)$
of $K= S(\Frac{\alpha}{\beta})$, $\Frac{\alpha}{\beta}>1$,
is obtained with $b=n+1$ where $n$ is the length of the continued fraction of
$\Frac\alpha\beta$ or $\Frac{\alpha}{\alpha-\beta}$ (see \cite{KP4}).
This allows us, using Algorithm \ref{cfs}, to know the minimal $b$ for which $\cC(3,b)$
is a projection of a given rational knot $K$.
\pn
In a similar manner, let $K= S(\Frac{\alpha}{\beta})$, $\Frac{\alpha}{\beta}>1$, $\beta$ even.
The minimal integer $b=2n+1$ for which there exists a continued fraction expansion
$r = [\eps_1, 2 \eps_2, \ldots, \eps_{2n-1},2\eps_{2n}]$, $\eps_i = \pm 1$,
such that $S(r)$ is equivalent to $K$,
is the smallest length of the continued fraction expansion $[\pm 1, \pm 2, \ldots, \pm 1, \pm2]$
of either $\Frac\alpha\beta$, $\Frac{\alpha}{2\alpha-\beta}$, $\Frac{\alpha}{\beta'}$ or
$\Frac{\alpha}{2\alpha-\beta'}$ where $0<\beta'<\alpha$, $\beta'$ even
and $\beta\beta'=\pm 1 \Mod \alpha$.
\pn
It happens that for some knots, there is a continued fraction expansion with smaller length
including $0$ instead of $\pm 2$.
The list of these knots (up to crossing number 10) is
$8_{12}, 9_{13}, 9_{15},9_{26},10_{8},10_{12},10_{13}, 10_{25}, 10_{29},10_{38},10_{42}$.
For example
for the knot $8_{12}= S(\Frac{29}{12})$, we have $\Frac{29}{12} = [1,0,1,2,1,0,1,2]$ while it
is not possible to get a shorter continued fraction corresponding to $8_{12}$.
We have enumerated all possible continued fraction expansions corresponding to diagrams
$\cC(4,b)$ to determine the minimal $b$ corresponding to a rational knot $K$.
\pn
In the next table we give the number of two-bridge knots of crossing number $N$
that have a projection $\cC(3,n)$ with $n\leq b$ and a projection $\cC(4,n')$ with $n'\leq b$.
\pn
\begin{center}
\setlength{\tabcolsep}{8pt}
\begin{tabular}{||c|cccccccc||cccccc||}
\hline
\multicolumn{15}{||c||}{\bf Minimal $b$}\\
\hline
&\multicolumn{8}{c||}{$a=3$}&\multicolumn{6}{c||}{$a=4$}\\
\hline
$b$ & 4 & 5 & 7 & 8 & 10 & 11 & 13 & 14&
3&5&7&9&11&13\\
\hline
$N=3$&1&1&1&1&1&1&1&1           &1&1&1&1&1&1\\
$N=4$& &1&1&1&1&1&1&1           & &1&1&1&1&1\\
$N=5$& & &2&2&2&2&2&2           & &2&2&2&2&2\\
$N=6$& & &1&3&3&3&3&3           & &1&3&3&3&3\\
$N=7$& & & &1&7&7&7&7           & & &5&7&7&7\\
$N=8$& & & & &3&12&12&12        & & &4&10&12&12\\
$N=9$& & & & &1&5&24&24         & & &1&14&24&24\\
$N=10$&& & & & &1&17&45         & & & &13&37&45\\
\hline
\end{tabular}
\end{center}
\pn
For example, looking at the 45 knots with crossing number 10 we see that: 1 of them
has a projection $\cC(3,11)$ (10 crossing points), 17 of them have a projection
$\cC(3,13)$  (12 crossing points) and
$45$ of them have a projection $\cC(3,14)$ (13 crossing points).
All of them have a projection $\cC(3,b)$ with $b \leq 14$.
We can also observe that for $b=10$ we have all knots with crossing number not greater than
6, 7 knots with crossing number 7, 3 knots with crossing number 8 and 1 knot with crossing number 9.
This last one is the Fibonacci knot $C(1,1,1,1,1,1,1,1,1) \sim \cC(3,10,17,0)$ (see \cite{KP3}).
\pn
Looking at the 45 knots with crossing number 10 we see that: 13 of them
have a projection $\cC(4,9)$ (12 crossing points), 37 of them have a projection
$\cC(4,11)$ (15 crossing points), 45 of them have a projection $\cC(4,13)$.
\pn
Let $K$ be a two-bridge knot. Once we know what we can expect as a diagram $\cC(a,b)$ for $K$,
we look for it as a Chebyshev knot $\cC(a,b,c,\phi)$. We see from the previous table that every
rational knot with crossing number $N\leq 10$ is a Chebyshev knot $\cC(3,b,c,r)$ with
$b \leq 14$ and a Chebyshev knot $\cC(4,b,c,r)$ with $b \leq 13$.
We proved in \cite{KP4} that a two-bridge knot with crossing number $N$ admits
a plane Chebyshev projection $\cC(3,b)$ with $b< \frac 32 N$.
In comparison, the number of crossing points for Lissajous
diagrams is far greater (see \cite{BDHZ}).
\pn
We have limited ourselves to the bounds $b\leq 21$, $c \leq 300$ and
$(b-1)(c-1) \leq 13 \cdot 299=3887$. The degrees of the
polynomials $R(\phi)$ giving the critical values are bounded by $3887$ when
$a=3$ and $5382$ when $a=4$.
A remarkable fact is that these polynomials have a large number of real roots
(in average 58\% when $a=3$ and $57\%$ when $a=4$, see Figure \ref{racs}, first column).
The proportion of non trivial knots $\cC(a,b,c,r)$ is approx. 25\% when $a=3$ and
39\% when $a=4$ (see column 2).
The proportion of non trivial distinct knots $\cC(a,b,c,\phi)$ is drawn in
Figure \ref{racs}, column 3.
\psfrag{x}{}\psfrag{y}{}
\pn
We conclude our paper by a list of the first
95 two-bridge knots (up to crossing number 10).
We give the Conway-Rolfsen numbering, their Schubert fraction (up to mirror symmetry) and
their presentation as Chebyshev knots. Most of them have a parametrization with the minimal
$b$. All of them have a parametrization with a minimal $(b-1)(c-1)$.

For example the knot $9_{5}$ admits the Chebyshev parametrization $\cC(3,17,45,1/364)$.
It is not minimal and we know that there is some other parametrization $\cC(3,13,c,r)$
where $c>300$ and $r$ is some rational number. This knot admits also the parametrization
$\cC(4,11,152,1/44)$ which is minimal with respect to $b$ and $\cC(4,19,22,1/20)$ that has minimal degree.

We get both $10_{20}$ and $10_{29}$ with minimal $b=14$ for the same value $c=292$. We had to compute
the polynomial $R_{3,14,292}$ of degree 3783 and 2185 real roots (see \ref{ex2}).
We obtain the knot $10_{20}=S(\frac{35}{11})$ with
$\cC(3,14,292,1/94)$ and the knot $10_{29}=S(\frac{63}{17})$ with $\cC(3,14,292,1/93)$.
\section{Conclusion}
We have shown that any two-bridge knot is a Chebyshev knot with $a=3$ and also with
$a=4$. For every $a,b,c$ integers ($a=3, 4$ and $a$, $b$ coprime), we have
described an algorithm that gives all Chebyshev knots $\cC(a,b,c,\phi)$.

Our experiments fully justify the use of certified algorithms and exact computations since
numerical methods would have certainly failed in finding for example  the knot
$10_{20}=S(\frac{35}{11})$ with $\cC(3,14,292,1/94)$ and the knot $10_{29}=S(\frac{63}{17})$
with $\cC(3,14,292,1/93)$. Also, an objective is now to consolidate and speed up our
algorithms in order to increase its capabilities.


As the zero-dimensional systems we study have a triangular structure,
we could try to get directly an exhaustive list of Chebyshev knots $\cC(a,b,c,\phi)$ without
computing additional resultants.
In case when $a$, $b$ and $c$ are relatively coprime, we can expect
that the real variety $V_{a,b,c}$ has only single points and try to get directly all possible signs.
In that case, our three black-boxes could be implemented using exclusively univariate functions that
compute recursively the roots of the systems to be solved without any additional rewriting.

\section*{Table}
Here is the list of the first 95 rational knots. We have given
Chebyshev parametrizations for $a=3$ and $a=4$. One corresponds to the minimal
$b$ and the other to the minimal total degree in $b,c$. For each parametrization
we give the corresponding Schubert fraction ($\alpha/\beta$),
the number of double points (DP) in the corresponding diagram
$\cC(a,b)$ so as the degree ($\deg$) of $V_{a,b,c}$. Note that sometimes we have fewer
double points with $a=4$ than with $a=3$. For 6 knots ($9_{5}$, $10_{36}$, $10_{39}$,
$10_{3}$, $10_{30}$, $10_{33}$), Chebyshev parametrizations
with $b$ minimal are not obtained with $c\leq 300$. Note that in \cite{KP4}, we have given
an algorithm that determines for any two-bridge knot $K$, the minimal integer $b$
and $C(t) \in \Q[t]$, $b+\deg C=3N$,
such that $x=T_3(t),\, y = T_b(t), z=C(t)$ is a parametrization of $K$.
\vfill
\pn
\hrule width 5cm height 2pt
\pn
Pierre-Vincent Koseleff, \\
INRIA-Paris-Rocquencourt Salsa \& Universit{\'e} Pierre et Marie Curie (UPMC-Paris 6)
\& Laboratoire d'Informatique de Paris 6, CNRS (UMR 7606)\\
e-mail: {\tt koseleff@math.jussieu.fr}
\pn
Daniel Pecker, \\
Universit{\'e} Pierre et Marie Curie (UPMC-Paris 6)\\
e-mail: {\tt pecker@math.jussieu.fr}
\pn
F. Rouillier, \\
INRIA-Paris-Rocquencourt Salsa \& Universit{\'e} Pierre et Marie Curie (UPMC-Paris 6)
\& Laboratoire d'Informatique de Paris 6, CNRS (UMR 7606)\\
e-mail: {\tt Fabrice.rouillier@inria.fr}
\setlength{\tabcolsep}{5pt}%
\renewcommand{\arraystretch}{.95}%
\small%
\def\bt{%
\eject%
\begin{center}%
\begin{tabular}{||c||c|r|c|c||c|r|c|c||}%
\hline%
\multicolumn{9}{||c||}{\bf Chebyshev parametrizations of the first rational knots}\rule[-7pt]{0pt}{20pt}\\%
\hline%
$K$&minimal $b$&$\alpha/\beta$&DP&$\deg$&min. $(b-1)(c-1)$&$\alpha/\beta$&DP&$\deg$\rule[-4pt]{0pt}{15pt}\\%
\hline}%
\def\et{\end{tabular}\end{center}}%
\bt
\multirow{2}{*}{$3_{1}$}&$\cC(3,4,5,0)$&$3/2$&$3$&$12$&$\cC(3,4,5,0)$&$3/2$&$3$&$12$\\ &$\cC(4,3,5,0)$&$3/2$&$3$&$12$&$\cC(4,5,6,1/5)$&$3/2$&$6$&$30$\\ \hline
\multirow{2}{*}{$4_{1}$}&$\cC(3,5,7,0)$&$5/3$&$4$&$24$&$\cC(3,5,7,0)$&$5/3$&$4$&$24$\\ &$\cC(4,5,12,1/23)$&$5/2$&$6$&$66$&$\cC(4,7,8,1/3)$&$5/2$&$9$&$63$\\ \hline
\multirow{2}{*}{$5_{1}$}&$\cC(3,7,8,0)$&$-5/4$&$6$&$42$&$\cC(3,7,8,0)$&$-5/4$&$6$&$42$\\ &$\cC(4,5,8,1/23)$&$-5/4$&$6$&$42$&$\cC(4,5,8,1/23)$&$-5/4$&$6$&$42$\\ \hline
\multirow{2}{*}{$5_{2}$}&$\cC(3,7,17,1/50)$&$-7/4$&$6$&$96$&$\cC(3,10,11,1/16)$&$7/4$&$9$&$90$\\ &$\cC(4,5,7,0)$&$-7/4$&$6$&$36$&$\cC(4,5,7,0)$&$-7/4$&$6$&$36$\\ \hline
\multirow{2}{*}{$6_{1}$}&$\cC(3,8,10,1/42)$&$-9/5$&$7$&$63$&$\cC(3,8,10,1/42)$&$-9/5$&$7$&$63$\\ &$\cC(4,7,16,1/39)$&$-9/4$&$9$&$135$&$\cC(4,7,16,1/39)$&$-9/4$&$9$&$135$\\ \hline
\multirow{2}{*}{$6_{2}$}&$\cC(3,8,19,1/46)$&$-11/7$&$7$&$126$&$\cC(3,8,19,1/46)$&$-11/7$&$7$&$126$\\ &$\cC(4,5,11,0)$&$11/8$&$6$&$60$&$\cC(4,5,11,0)$&$11/8$&$6$&$60$\\ \hline
\multirow{2}{*}{$6_{3}$}&$\cC(3,7,11,0)$&$13/8$&$6$&$60$&$\cC(3,7,11,0)$&$13/8$&$6$&$60$\\ &$\cC(4,7,36,1/42)$&$-13/8$&$9$&$315$&$\cC(4,9,14,1/29)$&$-13/8$&$12$&$156$\\ \hline
\multirow{2}{*}{$7_{1}$}&$\cC(3,10,11,0)$&$7/6$&$9$&$90$&$\cC(3,10,11,0)$&$7/6$&$9$&$90$\\ &$\cC(4,7,27,1/68)$&$7/6$&$9$&$234$&$\cC(4,9,12,1/18)$&$-7/6$&$12$&$132$\\ \hline
\multirow{2}{*}{$7_{2}$}&$\cC(3,10,27,1/50)$&$11/6$&$9$&$234$&$\cC(3,10,27,1/50)$&$11/6$&$9$&$234$\\ &$\cC(4,7,9,1/30)$&$11/6$&$9$&$72$&$\cC(4,7,9,1/30)$&$11/6$&$9$&$72$\\ \hline
\multirow{2}{*}{$7_{3}$}&$\cC(3,10,28,1/47)$&$13/4$&$9$&$243$&$\cC(3,10,28,1/47)$&$13/4$&$9$&$243$\\ &$\cC(4,7,27,1/80)$&$13/10$&$9$&$234$&$\cC(4,9,15,1/35)$&$13/10$&$12$&$168$\\ \hline
\multirow{2}{*}{$7_{4}$}&$\cC(3,10,36,1/306)$&$-15/4$&$9$&$315$&$\cC(3,11,27,1/238)$&$-15/11$&$10$&$260$\\ &$\cC(4,9,64,1/156)$&$-15/4$&$12$&$756$&$\cC(4,13,21,1/24)$&$-15/4$&$18$&$360$\\ \hline
\multirow{2}{*}{$7_{5}$}&$\cC(3,10,35,1/60)$&$17/12$&$9$&$306$&$\cC(3,13,14,1/24)$&$-17/10$&$12$&$156$\\ &$\cC(4,7,9,0)$&$17/10$&$9$&$72$&$\cC(4,7,9,0)$&$17/10$&$9$&$72$\\ \hline
\multirow{2}{*}{$7_{6}$}&$\cC(3,10,33,1/46)$&$19/8$&$9$&$288$&$\cC(3,14,15,1/26)$&$-19/7$&$13$&$182$\\ &$\cC(4,7,40,1/51)$&$19/8$&$9$&$351$&$\cC(4,9,13,5/44)$&$-19/8$&$12$&$144$\\ \hline
\multirow{2}{*}{$7_{7}$}&$\cC(3,8,13,0)$&$21/13$&$7$&$84$&$\cC(3,8,13,0)$&$21/13$&$7$&$84$\\ &$\cC(4,9,61,1/67)$&$-21/8$&$12$&$720$&$\cC(4,11,16,3/10)$&$-21/34$&$15$&$225$\\ \hline
\multirow{2}{*}{$8_{1}$}&$\cC(3,11,13,1/60)$&$13/7$&$10$&$120$&$\cC(3,11,13,1/60)$&$13/7$&$10$&$120$\\ &$\cC(4,9,118,1/67)$&$13/6$&$12$&$1404$&$\cC(4,15,20,2/15)$&$-13/24$&$21$&$399$\\ \hline
\multirow{2}{*}{$8_{2}$}&$\cC(3,11,28,1/80)$&$17/11$&$10$&$270$&$\cC(3,11,28,1/80)$&$17/11$&$10$&$270$\\ &$\cC(4,9,35,1/68)$&$17/14$&$12$&$408$&$\cC(4,13,17,1/6)$&$17/6$&$18$&$288$\\ \hline
\multirow{2}{*}{$8_{3}$}&$\cC(3,11,13,0)$&$17/13$&$10$&$120$&$\cC(3,11,13,0)$&$17/13$&$10$&$120$\\ &$\cC(4,11,101,1/85)$&$-17/4$&$15$&$1500$&$\cC(4,15,17,1/14)$&$-17/4$&$21$&$336$\\ \hline
\multirow{2}{*}{$8_{4}$}&$\cC(3,11,46,1/58)$&$19/15$&$10$&$450$&$\cC(3,17,22,1/144)$&$19/81$&$16$&$336$\\ &$\cC(4,7,25,1/75)$&$-19/14$&$9$&$216$&$\cC(4,11,12,3/14)$&$19/14$&$15$&$165$\\ \hline
\multirow{2}{*}{$8_{6}$}&$\cC(3,11,87,1/70)$&$23/13$&$10$&$860$&$\cC(3,17,22,1/44)$&$-23/13$&$16$&$336$\\ &$\cC(4,9,91,1/66)$&$23/16$&$12$&$1080$&$\cC(4,11,21,1/15)$&$-23/16$&$15$&$300$\\ \hline
\multirow{2}{*}{$8_{7}$}&$\cC(3,10,70,1/47)$&$-23/18$&$9$&$621$&$\cC(3,13,18,1/42)$&$23/18$&$12$&$204$\\ &$\cC(4,7,13,0)$&$-23/18$&$9$&$108$&$\cC(4,7,13,0)$&$-23/18$&$9$&$108$\\ \hline
\multirow{2}{*}{$8_{8}$}&$\cC(3,10,14,1/60)$&$-25/16$&$9$&$117$&$\cC(3,10,14,1/60)$&$-25/16$&$9$&$117$\\ &$\cC(4,7,24,1/236)$&$-25/14$&$9$&$207$&$\cC(4,7,24,1/236)$&$-25/14$&$9$&$207$\\ \hline
\multirow{2}{*}{$8_{9}$}&$\cC(3,11,110,1/168)$&$25/7$&$10$&$1090$&$\cC(3,17,25,1/96)$&$-25/57$&$16$&$384$\\ &$\cC(4,7,16,1/54)$&$-25/18$&$9$&$135$&$\cC(4,7,16,1/54)$&$-25/18$&$9$&$135$\\ \hline
\multirow{2}{*}{$8_{11}$}&$\cC(3,11,100,1/84)$&$27/17$&$10$&$990$&$\cC(3,14,38,1/66)$&$-27/17$&$13$&$481$\\ &$\cC(4,9,30,1/82)$&$27/8$&$12$&$348$&$\cC(4,11,13,2/25)$&$-27/8$&$15$&$180$\\ \hline
\multirow{2}{*}{$8_{12}$}&$\cC(3,11,54,1/152)$&$29/17$&$10$&$530$&$\cC(3,14,41,1/19)$&$29/17$&$13$&$520$\\ &$\cC(4,9,103,1/68)$&$29/12$&$12$&$1224$&$\cC(4,15,18,4/29)$&$-29/12$&$21$&$357$\\ \hline
\et
\bt
\multirow{2}{*}{$8_{13}$}&$\cC(3,10,17,1/62)$&$-29/8$&$9$&$144$&$\cC(3,10,17,1/62)$&$-29/8$&$9$&$144$\\ &$\cC(4,9,104,1/66)$&$29/18$&$12$&$1236$&$\cC(4,13,18,17/45)$&$-29/8$&$18$&$306$\\ \hline
\multirow{2}{*}{$8_{14}$}&$\cC(3,11,93,1/86)$&$31/13$&$10$&$920$&$\cC(3,17,22,1/26)$&$-31/19$&$16$&$336$\\ &$\cC(4,9,47,1/200)$&$-31/44$&$12$&$552$&$\cC(4,13,18,1/20)$&$-31/18$&$18$&$306$\\ \hline
\multirow{2}{*}{$9_{1}$}&$\cC(3,13,14,0)$&$-9/8$&$12$&$156$&$\cC(3,13,14,0)$&$-9/8$&$12$&$156$\\ &$\cC(4,9,87,1/59)$&$-9/8$&$12$&$1032$&$\cC(4,13,14,3/43)$&$9/8$&$18$&$234$\\ \hline
\multirow{2}{*}{$9_{2}$}&$\cC(3,13,37,1/114)$&$-15/2$&$12$&$432$&$\cC(3,19,24,3/82)$&$15/8$&$18$&$414$\\ &$\cC(4,11,33,1/15)$&$-15/8$&$15$&$480$&$\cC(4,17,20,7/23)$&$-15/2$&$24$&$456$\\ \hline
\multirow{2}{*}{$9_{3}$}&$\cC(3,13,143,1/98)$&$-19/6$&$12$&$1704$&$\cC(3,16,46,1/54)$&$19/6$&$15$&$675$\\ &$\cC(4,9,159,1/108)$&$-19/16$&$12$&$1896$&$\cC(4,17,18,1/23)$&$19/44$&$24$&$408$\\ \hline
\multirow{2}{*}{$9_{4}$}&$\cC(3,13,115,1/164)$&$-21/4$&$12$&$1368$&$\cC(3,19,24,1/40)$&$21/16$&$18$&$414$\\ &$\cC(4,9,106,3/301)$&$-21/16$&$12$&$1260$&$\cC(4,17,18,3/17)$&$-21/4$&$24$&$408$\\ \hline
\multirow{2}{*}{$9_{5}$}&$b=13$&$23/4$&$12$&$12c-12$&$\cC(3,17,45,1/182)$&$23/29$&$16$&$704$\\ &$\cC(4,11,152,1/44)$&$23/6$&$15$&$2265$&$\cC(4,19,22,1/20)$&$23/6$&$27$&$567$\\ \hline
\multirow{2}{*}{$9_{6}$}&$\cC(3,13,64,1/102)$&$-27/22$&$12$&$756$&$\cC(3,16,17,1/34)$&$27/16$&$15$&$240$\\ &$\cC(4,9,11,1/55)$&$-27/16$&$12$&$120$&$\cC(4,9,11,1/55)$&$-27/16$&$12$&$120$\\ \hline
\multirow{2}{*}{$9_{7}$}&$\cC(3,13,116,1/80)$&$-29/20$&$12$&$1380$&$\cC(3,22,24,3/26)$&$-29/16$&$21$&$483$\\ &$\cC(4,9,201,1/131)$&$-29/16$&$12$&$2400$&$\cC(4,11,17,1/65)$&$-29/16$&$15$&$240$\\ \hline
\multirow{2}{*}{$9_{8}$}&$\cC(3,13,121,1/170)$&$-31/20$&$12$&$1440$&$\cC(3,16,25,3/97)$&$31/14$&$15$&$360$\\ &$\cC(4,11,187,1/147)$&$31/14$&$15$&$2790$&$\cC(4,17,20,6/19)$&$-31/14$&$24$&$456$\\ \hline
\multirow{2}{*}{$9_{9}$}&$\cC(3,13,123,1/80)$&$-31/22$&$12$&$1464$&$\cC(3,19,52,2/29)$&$-31/24$&$18$&$918$\\ &$\cC(4,9,31,1/66)$&$-31/22$&$12$&$360$&$\cC(4,9,31,1/66)$&$-31/22$&$12$&$360$\\ \hline
\multirow{2}{*}{$9_{10}$}&$\cC(3,13,246,1/110)$&$-33/10$&$12$&$2940$&$\cC(3,19,53,1/32)$&$33/10$&$18$&$936$\\ &$\cC(4,11,29,1/45)$&$33/10$&$15$&$420$&$\cC(4,11,29,1/45)$&$33/10$&$15$&$420$\\ \hline
\multirow{2}{*}{$9_{11}$}&$\cC(3,13,114,1/106)$&$-33/26$&$12$&$1356$&$\cC(3,23,24,1/10)$&$-33/19$&$22$&$506$\\ &$\cC(4,9,33,1/75)$&$-33/26$&$12$&$384$&$\cC(4,9,33,1/75)$&$-33/26$&$12$&$384$\\ \hline
\multirow{2}{*}{$9_{12}$}&$\cC(3,13,36,1/80)$&$-35/22$&$12$&$420$&$\cC(3,13,36,1/80)$&$-35/22$&$12$&$420$\\ &$\cC(4,11,68,1/77)$&$35/22$&$15$&$1005$&$\cC(4,15,24,1/17)$&$-35/22$&$21$&$483$\\ \hline
\multirow{2}{*}{$9_{13}$}&$\cC(3,13,53,1/78)$&$-37/26$&$12$&$624$&$\cC(3,16,31,1/42)$&$37/10$&$15$&$450$\\ &$\cC(4,9,41,1/91)$&$-37/26$&$12$&$480$&$\cC(4,17,18,1/6)$&$-37/26$&$24$&$408$\\ \hline
\multirow{2}{*}{$9_{14}$}&$\cC(3,11,83,1/74)$&$-37/23$&$10$&$820$&$\cC(3,17,18,1/9)$&$-37/29$&$16$&$272$\\ &$\cC(4,11,176,1/108)$&$37/14$&$15$&$2625$&$\cC(4,15,22,4/11)$&$37/8$&$21$&$441$\\ \hline
\multirow{2}{*}{$9_{15}$}&$\cC(3,13,144,1/310)$&$39/22$&$12$&$1716$&$\cC(3,17,18,1/34)$&$39/17$&$16$&$272$\\ &$\cC(4,9,39,1/85)$&$-39/22$&$12$&$456$&$\cC(4,13,15,7/45)$&$39/16$&$18$&$252$\\ \hline
\multirow{2}{*}{$9_{17}$}&$\cC(3,11,16,0)$&$-39/25$&$10$&$150$&$\cC(3,11,16,0)$&$-39/25$&$10$&$150$\\ &$\cC(4,9,92,1/92)$&$-39/14$&$12$&$1092$&$\cC(4,11,24,1/53)$&$39/14$&$15$&$345$\\ \hline
\multirow{2}{*}{$9_{18}$}&$\cC(3,13,194,1/144)$&$-41/12$&$12$&$2316$&$\cC(3,16,43,1/36)$&$-41/24$&$15$&$630$\\ &$\cC(4,9,11,0)$&$-41/24$&$12$&$120$&$\cC(4,9,11,0)$&$-41/24$&$12$&$120$\\ \hline
\multirow{2}{*}{$9_{19}$}&$\cC(3,11,83,1/82)$&$-41/25$&$10$&$820$&$\cC(3,14,22,1/110)$&$41/105$&$13$&$273$\\ &$\cC(4,9,32,1/416)$&$41/64$&$12$&$372$&$\cC(4,9,32,1/416)$&$41/64$&$12$&$372$\\ \hline
\multirow{2}{*}{$9_{20}$}&$\cC(3,13,275,1/86)$&$-41/26$&$12$&$3288$&$\cC(3,16,46,1/34)$&$41/26$&$15$&$675$\\ &$\cC(4,7,17,0)$&$41/30$&$9$&$144$&$\cC(4,7,17,0)$&$41/30$&$9$&$144$\\ \hline
\multirow{2}{*}{$9_{21}$}&$\cC(3,13,179,1/106)$&$-43/18$&$12$&$2136$&$\cC(3,16,29,9/172)$&$-43/18$&$15$&$420$\\ &$\cC(4,9,94,1/105)$&$43/68$&$12$&$1116$&$\cC(4,11,23,1/23)$&$43/18$&$15$&$330$\\ \hline
\multirow{2}{*}{$9_{23}$}&$\cC(3,13,44,1/98)$&$-45/26$&$12$&$516$&$\cC(3,13,44,1/98)$&$-45/26$&$12$&$516$\\ &$\cC(4,13,16,1/20)$&$-45/26$&$18$&$270$&$\cC(4,13,16,1/20)$&$-45/26$&$18$&$270$\\ \hline
\multirow{2}{*}{$9_{26}$}&$\cC(3,11,25,1/92)$&$-47/29$&$10$&$240$&$\cC(3,11,25,1/92)$&$-47/29$&$10$&$240$\\ &$\cC(4,9,184,1/79)$&$-47/34$&$12$&$2196$&$\cC(4,13,17,1/25)$&$47/18$&$18$&$288$\\ \hline
\et
\bt
\multirow{2}{*}{$9_{27}$}&$\cC(3,13,180,1/84)$&$-49/18$&$12$&$2148$&$\cC(3,17,51,1/30)$&$49/31$&$16$&$800$\\ &$\cC(4,9,39,1/66)$&$-49/30$&$12$&$456$&$\cC(4,15,20,1/28)$&$49/80$&$21$&$399$\\ \hline
\multirow{2}{*}{$9_{31}$}&$\cC(3,10,17,0)$&$55/34$&$9$&$144$&$\cC(3,10,17,0)$&$55/34$&$9$&$144$\\ &$\cC(4,11,68,1/100)$&$55/34$&$15$&$1005$&$\cC(4,13,24,1/43)$&$-55/34$&$18$&$414$\\ \hline
\multirow{2}{*}{$10_{1}$}&$\cC(3,14,38,1/102)$&$-17/9$&$13$&$481$&$\cC(3,17,19,1/44)$&$17/9$&$16$&$288$\\ &$\cC(4,11,141,1/44)$&$-17/8$&$15$&$2100$&$\cC(4,19,23,1/33)$&$17/8$&$27$&$594$\\ \hline
\multirow{2}{*}{$10_{2}$}&$\cC(3,14,43,1/176)$&$-23/3$&$13$&$546$&$\cC(3,20,25,9/314)$&$23/15$&$19$&$456$\\ &$\cC(4,9,278,1/85)$&$23/20$&$12$&$3324$&$\cC(4,17,20,4/61)$&$-23/38$&$24$&$456$\\ \hline
\multirow{2}{*}{$10_{3}$}&$\cC(3,14,16,1/94)$&$-25/19$&$13$&$195$&$\cC(3,14,16,1/94)$&$-25/19$&$13$&$195$\\ &$b=13$&$25/4$&$18$&$18c-18$&$\cC(4,19,20,2/17)$&$25/6$&$27$&$513$\\ \hline
\multirow{2}{*}{$10_{4}$}&$\cC(3,14,101,1/130)$&$-27/23$&$13$&$1300$&$\cC(3,23,26,1/38)$&$27/23$&$22$&$550$\\ &$\cC(4,9,257,1/145)$&$27/20$&$12$&$3072$&$\cC(4,11,34,2/17)$&$27/20$&$15$&$495$\\ \hline
\multirow{2}{*}{$10_{5}$}&$\cC(3,13,169,9/1034)$&$33/28$&$12$&$2016$&$\cC(3,19,24,1/116)$&$33/28$&$18$&$414$\\ &$\cC(4,11,194,1/110)$&$-33/28$&$15$&$2895$&$\cC(4,13,20,2/35)$&$33/28$&$18$&$342$\\ \hline
\multirow{2}{*}{$10_{6}$}&$\cC(3,14,128,1/92)$&$-37/7$&$13$&$1651$&$\cC(3,20,25,1/42)$&$37/21$&$19$&$456$\\ &$\cC(4,13,72,1/47)$&$37/16$&$18$&$1278$&$\cC(4,17,20,1/15)$&$-37/58$&$24$&$456$\\ \hline
\multirow{2}{*}{$10_{7}$}&$\cC(3,14,127,1/128)$&$-43/27$&$13$&$1638$&$\cC(3,23,33,1/24)$&$43/27$&$22$&$704$\\ &$\cC(4,11,229,1/70)$&$-43/16$&$15$&$3420$&$\cC(4,19,21,1/17)$&$-43/8$&$27$&$540$\\ \hline
\multirow{2}{*}{$10_{8}$}&$\cC(3,14,37,1/144)$&$-29/23$&$13$&$468$&$\cC(3,14,37,1/144)$&$-29/23$&$13$&$468$\\ &$\cC(4,11,77,1/66)$&$-29/24$&$15$&$1140$&$\cC(4,13,36,2/45)$&$29/24$&$18$&$630$\\ \hline
\multirow{2}{*}{$10_{9}$}&$\cC(3,14,281,1/232)$&$-39/11$&$13$&$3640$&$\cC(3,17,43,1/186)$&$39/89$&$16$&$672$\\ &$\cC(4,9,35,1/133)$&$39/28$&$12$&$408$&$\cC(4,13,17,17/120)$&$39/28$&$18$&$288$\\ \hline
\multirow{2}{*}{$10_{10}$}&$\cC(3,13,253,1/250)$&$45/28$&$12$&$3024$&$\cC(3,19,31,1/31)$&$-45/118$&$18$&$540$\\ &$\cC(4,11,102,1/45)$&$-45/28$&$15$&$1515$&$\cC(4,15,20,1/10)$&$45/28$&$21$&$399$\\ \hline
\multirow{2}{*}{$10_{11}$}&$\cC(3,14,101,1/116)$&$-43/33$&$13$&$1300$&$\cC(3,20,25,1/196)$&$-43/185$&$19$&$456$\\ &$\cC(4,11,126,1/97)$&$-43/30$&$15$&$1875$&$\cC(4,21,24,2/31)$&$-43/76$&$30$&$690$\\ \hline
\multirow{2}{*}{$10_{12}$}&$\cC(3,13,61,1/178)$&$47/36$&$12$&$720$&$\cC(3,19,24,1/44)$&$47/36$&$18$&$414$\\ &$\cC(4,11,115,1/393)$&$-47/58$&$15$&$1710$&$\cC(4,13,20,1/18)$&$47/36$&$18$&$342$\\ \hline
\multirow{2}{*}{$10_{13}$}&$\cC(3,14,211,3/322)$&$-53/41$&$13$&$2730$&$\cC(3,23,26,1/26)$&$53/31$&$22$&$550$\\ &$\cC(4,11,147,1/84)$&$-53/22$&$15$&$2190$&$\cC(4,15,18,3/22)$&$-53/22$&$21$&$357$\\ \hline
\multirow{2}{*}{$10_{14}$}&$\cC(3,14,139,1/180)$&$-57/13$&$13$&$1794$&$\cC(3,17,44,2/43)$&$-57/13$&$16$&$688$\\ &$\cC(4,9,35,2/135)$&$57/44$&$12$&$408$&$\cC(4,9,35,2/135)$&$57/44$&$12$&$408$\\ \hline
\multirow{2}{*}{$10_{15}$}&$\cC(3,13,17,1/80)$&$43/24$&$12$&$192$&$\cC(3,13,17,1/80)$&$43/24$&$12$&$192$\\ &$\cC(4,9,105,1/40)$&$43/24$&$12$&$1248$&$\cC(4,17,27,2/35)$&$43/24$&$24$&$624$\\ \hline
\multirow{2}{*}{$10_{16}$}&$\cC(3,14,127,1/104)$&$-47/37$&$13$&$1638$&$\cC(3,23,25,1/9)$&$-47/37$&$22$&$528$\\ &$\cC(4,11,37,1/57)$&$-47/14$&$15$&$540$&$\cC(4,19,21,1/20)$&$-47/14$&$27$&$540$\\ \hline
\multirow{2}{*}{$10_{17}$}&$\cC(3,13,194,1/79)$&$41/32$&$12$&$2316$&$\cC(3,19,24,1/220)$&$41/32$&$18$&$414$\\ &$\cC(4,9,16,1/69)$&$41/32$&$12$&$180$&$\cC(4,9,16,1/69)$&$41/32$&$12$&$180$\\ \hline
\multirow{2}{*}{$10_{18}$}&$\cC(3,14,37,1/148)$&$-55/43$&$13$&$468$&$\cC(3,14,37,1/148)$&$-55/43$&$13$&$468$\\ &$\cC(4,11,211,1/52)$&$-55/32$&$15$&$3150$&$\cC(4,19,21,2/11)$&$55/32$&$27$&$540$\\ \hline
\multirow{2}{*}{$10_{19}$}&$\cC(3,13,128,1/158)$&$51/14$&$12$&$1524$&$\cC(3,22,33,1/48)$&$51/142$&$21$&$672$\\ &$\cC(4,9,162,1/181)$&$51/40$&$12$&$1932$&$\cC(4,13,23,1/148)$&$-51/40$&$18$&$396$\\ \hline
\multirow{2}{*}{$10_{20}$}&$\cC(3,14,292,1/94)$&$-35/11$&$13$&$3783$&$\cC(3,23,29,1/22)$&$-35/19$&$22$&$616$\\ &$\cC(4,11,298,1/133)$&$-35/16$&$15$&$4455$&$\cC(4,19,28,1/12)$&$35/54$&$27$&$729$\\ \hline
\multirow{2}{*}{$10_{21}$}&$\cC(3,14,133,1/108)$&$-45/31$&$13$&$1716$&$\cC(3,20,46,5/72)$&$-45/29$&$19$&$855$\\ &$\cC(4,11,193,1/60)$&$-45/16$&$15$&$2880$&$\cC(4,19,23,1/36)$&$45/16$&$27$&$594$\\ \hline
\multirow{2}{*}{$10_{22}$}&$\cC(3,14,230,1/554)$&$49/15$&$13$&$2977$&$\cC(3,20,27,1/126)$&$-49/15$&$19$&$494$\\ &$\cC(4,9,96,1/52)$&$49/36$&$12$&$1140$&$\cC(4,11,25,2/31)$&$49/36$&$15$&$360$\\ \hline
\et
\bt
\multirow{2}{*}{$10_{23}$}&$\cC(3,13,124,1/362)$&$-59/18$&$12$&$1476$&$\cC(3,19,26,1/64)$&$59/18$&$18$&$450$\\ &$\cC(4,11,38,1/105)$&$-59/18$&$15$&$555$&$\cC(4,19,20,2/13)$&$59/36$&$27$&$513$\\ \hline
\multirow{2}{*}{$10_{24}$}&$\cC(3,14,127,3/319)$&$-55/31$&$13$&$1638$&$\cC(3,23,36,5/82)$&$55/31$&$22$&$770$\\ &$\cC(4,11,247,1/74)$&$-55/24$&$15$&$3690$&$\cC(4,19,23,1/37)$&$55/24$&$27$&$594$\\ \hline
\multirow{2}{*}{$10_{25}$}&$\cC(3,14,148,1/108)$&$-65/41$&$13$&$1911$&$\cC(3,17,64,1/46)$&$-65/19$&$16$&$1008$\\ &$\cC(4,9,116,2/135)$&$65/46$&$12$&$1380$&$\cC(4,15,20,1/43)$&$65/106$&$21$&$399$\\ \hline
\multirow{2}{*}{$10_{26}$}&$\cC(3,14,110,1/98)$&$-61/17$&$13$&$1417$&$\cC(3,17,67,1/86)$&$-61/17$&$16$&$1056$\\ &$\cC(4,9,35,1/77)$&$61/44$&$12$&$408$&$\cC(4,9,35,1/77)$&$61/44$&$12$&$408$\\ \hline
\multirow{2}{*}{$10_{27}$}&$\cC(3,13,126,1/218)$&$71/50$&$12$&$1500$&$\cC(3,22,36,1/50)$&$71/44$&$21$&$735$\\ &$\cC(4,11,278,1/115)$&$-71/44$&$15$&$4155$&$\cC(4,15,26,2/71)$&$71/44$&$21$&$525$\\ \hline
\multirow{2}{*}{$10_{28}$}&$\cC(3,13,191,1/112)$&$53/14$&$12$&$2280$&$\cC(3,19,25,5/138)$&$53/34$&$18$&$432$\\ &$\cC(4,11,114,1/139)$&$53/92$&$15$&$1695$&$\cC(4,15,38,1/38)$&$53/92$&$21$&$777$\\ \hline
\multirow{2}{*}{$10_{29}$}&$\cC(3,14,292,1/93)$&$-63/17$&$13$&$3783$&$\cC(3,17,45,1/42)$&$-63/37$&$16$&$704$\\ &$\cC(4,9,168,1/106)$&$63/46$&$12$&$2004$&$\cC(4,15,23,1/8)$&$63/46$&$21$&$462$\\ \hline
\multirow{2}{*}{$10_{30}$}&$\cC(3,14,201,1/96)$&$-67/41$&$13$&$2600$&$\cC(3,17,39,1/34)$&$-67/49$&$16$&$608$\\ &$b=13$&$67/18$&$18$&$18c-18$&$\cC(4,21,25,3/77)$&$67/18$&$30$&$720$\\ \hline
\multirow{2}{*}{$10_{31}$}&$\cC(3,13,103,1/80)$&$57/32$&$12$&$1224$&$\cC(3,19,27,1/32)$&$-57/32$&$18$&$468$\\ &$\cC(4,9,111,1/66)$&$57/32$&$12$&$1320$&$\cC(4,15,20,1/11)$&$57/32$&$21$&$399$\\ \hline
\multirow{2}{*}{$10_{32}$}&$\cC(3,14,148,1/172)$&$-69/19$&$13$&$1911$&$\cC(3,16,56,1/166)$&$-69/50$&$15$&$825$\\ &$\cC(4,11,134,1/103)$&$-69/50$&$15$&$1995$&$\cC(4,15,22,1/34)$&$-69/40$&$21$&$441$\\ \hline
\multirow{2}{*}{$10_{33}$}&$\cC(3,13,182,1/105)$&$65/18$&$12$&$2172$&$\cC(3,22,40,1/38)$&$65/148$&$21$&$819$\\ &$b=13$&$65/18$&$18$&$18c-18$&$\cC(4,25,30,5/17)$&$65/18$&$36$&$1044$\\ \hline
\multirow{2}{*}{$10_{34}$}&$\cC(3,13,41,1/90)$&$37/20$&$12$&$480$&$\cC(3,16,20,1/44)$&$-37/24$&$15$&$285$\\ &$\cC(4,11,142,1/122)$&$-37/20$&$15$&$2115$&$\cC(4,13,19,1/11)$&$-37/20$&$18$&$324$\\ \hline
\multirow{2}{*}{$10_{35}$}&$\cC(3,14,38,1/108)$&$-49/27$&$13$&$481$&$\cC(3,14,38,1/108)$&$-49/27$&$13$&$481$\\ &$\cC(4,13,273,3/697)$&$49/20$&$18$&$4896$&$\cC(4,19,24,1/9)$&$-49/78$&$27$&$621$\\ \hline
\multirow{2}{*}{$10_{36}$}&$b=14$&$51/20$&$13$&$13c-13$&$\cC(3,20,32,1/44)$&$51/31$&$19$&$589$\\ &$\cC(4,9,179,1/222)$&$51/28$&$12$&$2136$&$\cC(4,11,17,1/77)$&$-51/28$&$15$&$240$\\ \hline
\multirow{2}{*}{$10_{37}$}&$\cC(3,13,17,0)$&$53/30$&$12$&$192$&$\cC(3,13,17,0)$&$53/30$&$12$&$192$\\ &$\cC(4,13,19,5/59)$&$-53/30$&$18$&$324$&$\cC(4,13,19,5/59)$&$-53/30$&$18$&$324$\\ \hline
\multirow{2}{*}{$10_{38}$}&$\cC(3,14,120,1/134)$&$-59/25$&$13$&$1547$&$\cC(3,17,23,5/122)$&$-59/33$&$16$&$352$\\ &$\cC(4,13,72,1/43)$&$59/26$&$18$&$1278$&$\cC(4,23,27,2/19)$&$-59/34$&$33$&$858$\\ \hline
\multirow{2}{*}{$10_{39}$}&$b=14$&$61/22$&$13$&$13c-13$&$\cC(3,20,25,1/34)$&$61/39$&$19$&$456$\\ &$\cC(4,9,277,1/119)$&$61/36$&$12$&$3312$&$\cC(4,11,39,1/18)$&$-61/36$&$15$&$570$\\ \hline
\multirow{2}{*}{$10_{40}$}&$\cC(3,13,190,1/112)$&$75/44$&$12$&$2268$&$\cC(3,16,34,1/21)$&$75/46$&$15$&$495$\\ &$\cC(4,11,57,1/291)$&$75/106$&$15$&$840$&$\cC(4,17,23,2/25)$&$-75/44$&$24$&$528$\\ \hline
\multirow{2}{*}{$10_{41}$}&$\cC(3,14,208,1/110)$&$-71/41$&$13$&$2691$&$\cC(3,17,55,1/23)$&$-71/41$&$16$&$864$\\ &$\cC(4,9,165,2/141)$&$-71/112$&$12$&$1968$&$\cC(4,15,25,6/47)$&$71/30$&$21$&$504$\\ \hline
\multirow{2}{*}{$10_{42}$}&$\cC(3,13,134,1/166)$&$81/50$&$12$&$1596$&$\cC(3,16,25,1/22)$&$81/50$&$15$&$360$\\ &$\cC(4,11,131,1/132)$&$-81/34$&$15$&$1950$&$\cC(4,13,32,5/47)$&$81/34$&$18$&$558$\\ \hline
\multirow{2}{*}{$10_{43}$}&$\cC(3,13,174,1/114)$&$73/46$&$12$&$2076$&$\cC(3,16,51,15/976)$&$-73/46$&$15$&$750$\\ &$\cC(4,11,277,1/134)$&$-73/46$&$15$&$4140$&$\cC(4,13,36,1/48)$&$73/100$&$18$&$630$\\ \hline
\multirow{2}{*}{$10_{44}$}&$\cC(3,14,132,1/114)$&$-79/49$&$13$&$1703$&$\cC(3,20,53,5/148)$&$79/29$&$19$&$988$\\ &$\cC(4,11,64,1/119)$&$-79/30$&$15$&$945$&$\cC(4,15,24,3/13)$&$-79/50$&$21$&$483$\\ \hline
\multirow{2}{*}{$10_{45}$}&$\cC(3,11,19,0)$&$89/55$&$10$&$180$&$\cC(3,11,19,0)$&$89/55$&$10$&$180$\\ &$\cC(4,13,132,1/238)$&$89/34$&$18$&$2358$&$\cC(4,15,26,2/21)$&$89/144$&$21$&$525$\\ \hline
\et
\end{document}